\documentclass[a4paper,11pt]{article}
\makeatletter
   \let\accent@spacefactor\relax
 \makeatother       
\usepackage[french]{babel}
\usepackage{amsmath,amscd,amssymb}

\def \cC {\mathfrak{C}}
\def \cv {{\check V}}
\def \b {\beta}
\def \Z {{\mathbb Z}}
\def \bC {{\mathbb C}}
\def \F {{\bf F}}
\def \Q {{\bf Q}}
\def \q {{\bf q}}
\def \aaq {\overline{\q}}
\def \LV {\Lambda^2 V}
\def \p {{\mathbb P}}
\def \pu {{\mathbb P}^1}
\def \opu {\oo_{\pu}}
\def \noi {\noindent}
\def \vs {\vskip}
\def \I {{\cal I}}
\def \oo {{\cal O}}
\def \P {{\bf \Psi}}
\def \G {{\mathbb G}}
\def \MD {{\cal M}(d)}

\def \Mor #1#2{{\rm{Mor}}_{#1}(\pu,#2)}
\def \fl {\longrightarrow}

\def \R {{\bf R}}

\def \hom #1#2{{\rm{Hom}}(#1,#2)}
\def \a {\alpha}

\begin{document}

\newtheorem{theo}{Th\'eor\`eme}
\newtheorem{exe}{Exemple}
\newtheorem{pr}{Proposition}
\newtheorem{defi}{D\'efinition}
\newtheorem{prob}{Probl\`eme}
\newtheorem{rem}{Remarque}
\newtheorem{lem}{Lemme}
\newtheorem{cor}{Corollaire}

\renewcommand{\thesection}{\arabic{section}}

\centerline{\Large{\bf Lieu singulier des surfaces rationnelles
 r\'egl\'ees}}

\vs 0.1 cm

\centerline{\large{\textsc{Nicolas Perrin}}}
\centerline{{Universit\'e de Versailles}}
\vs -0.1 cm
\centerline{45 avenue des \'Etats-Unis}
\vs -0.1 cm
\centerline{78035 Versailles Cedex}
\vs -0.1 cm
\centerline{email : \texttt{perrin@math.uvsq.fr}}

\vs 0.5 cm

\def \pv {{\mathbb P}(V)}
\def \E {{\cal E}}
\def \Gr {{\mathbb G}r}
\def \V {{\cal V}}
\def \opv {{\oo_{\pv}}}
\def \S {{S}}
\def \ps {\p(S_2)}
\def \ops {{\oo_{\ps}}}
\def \e {{\bf Ext}}
\def \h {{\bf Hom}}
\def \dm {{\textit{D\'emonstration}}}
\def \prop {{\bf Proposition}}
\def \g {{\mathbb G}}
\def \H {{\bf H}}
\def \L {{\Lambda}}
\def \oopu {{\otimes\opu}}
\def \oopv {{\otimes\opv}}
\def \oops {{\otimes\ops}}
\def \pd {\p({\S_d})}
\def \pvd {\p({\check V})}
\def \plv {\p(\LV)}
\def \s {{s^2}}
\def \ot {{\otimes}}
\def \w {{\wedge}}
\def \C {{\cal C}}
\def \cl {{\cal L}}
\def \pru {{\rm{pr}_1}}
\def \prd {{\rm{pr}_2}}
\def \rdp {{\bf{\cal R}_d}(\pv)}
\def \rdg {{\bf{\cal R}_d}(\G)}
\def \hdv {{\bf H}_d(V)}
\def \M {{\bf Hom}}
\def \Md {\Mor{d}{\G}}
\def \vp {\varphi}
\def \aq {{{\overline {\bf Q}}}}
\def \crp {{\mathfrak P}}
\def \psif {\psi_{{\rm{f}}}}
\def \d {[\frac{d}{2}]}
\def \Som {S_a\oplus S_{d-a}}
\def \som {\opu(a)\oplus\opu(d-a)}
\def \pt {\widetilde{\psi}}

\vs 0.3 cm

\centerline{\begin{minipage}{13 cm}
\footnotesize{\tableofcontents}
\end{minipage}}


\section*{Introduction}
\addcontentsline{toc}{section}{Introduction}

Dans cet article, nous appelons surface r\'egl\'ee toute surface de
$\p^3$ recouverte par une famille de dimension $1$ de droites de
l'espace, ce sont les g\'en\'eratrices de la surface. Cette
d\'efinition correpond \`a la notion
historique de surface r\'egl\'ee. Une telle surface est d\'efinie par
une courbe trac\'ee sur la grassmannienne $\G$ des droites de
$\p^3$. La surface est rationnelle si la courbe de $\G$ l'est. Le lieu
singulier de la surface est l'ensemble des points d'intersection de
deux g\'en\'eratrices.

Dans \cite{P1}, nous avons constat\'e que dans le cas des surfaces
quintiques rationnelles r\'egl\'ees, le lieu singulier de la surface
d\'etermine compl\`etement la surface. Nous nous proposons dans cet
article d'\'etudier le probl\`eme suivant :

\vs 0.1 cm

\noi
{\bf Probl\`eme} : \textit{Toute surface rationelle r\'egl\'ee
est-elle d\'etermin\'ee par son lieu singulier ?}

\vs 0.1 cm

Nous \'etudierons ce probl\`eme pour les surface rationnelles
r\'egl\'ees \textit{param\'etr\'ees}. L'ensemble de ces surfaces est
repr\'esent\'e par le sch\'ema $\Mor{d}{\G}$ des morphismes de degr\'e
$d$ de $\pu$ dans $\G$ (voir \cite{GR}). Pour tout morphisme $f$, la
surface correspondante est d\'efinie par la courbe $f(\pu)$ et le
morphisme nous donne un param\'etrage de cette courbe par $\pu$. La
vari\'et\'e $\Mor{d}{\G}$
est irr\'eductible et lisse de dimension $4d+4$ (voir \cite{P2}). Elle
est munie d'actions de $PGL_2$ et de $PO(q)$ o\`u $q$ est la forme
quadratique d\'efinissant $\G$ dans $\p^5$. Nous identifierons $PGL_4$
\`a $PSO(q)$. Le lieu singulier \textit{abstrait} d'une surface rationnelle
param\'etr\'ee $f\in\Md$ est d\'efini par :

\vs 0.1 cm

\noi
{\bf D\'efinition} : \textit{Le lieu singulier abstrait de la surface $f$ est l'ensemble des paires
de g\'en\'era\-trices qui se coupent, on le note $\P(f)$.}

\vs 0.1 cm

C'est en g\'en\'eral une courbe de degr\'e $d-2$ de $S^2\pu$. Une
surface rationelle r\'egl\'ee de degr\'e inf\'erieur \`a $2$
g\'en\'erale est lisse. Nous supposerons donc que $d\geq 3$. De plus, si la
surface $f$ ou sa duale est un c\^one, alors son lieu singulier
abstrait est $S^2\pu$ tout entier. Notons $\R_d$ l'ouvert de
$\Mor{d}{\G}$ des morphismes $f$ tels que l'image $f(\pu)$ n'est pas
contenue dans un plan totalement isotrope pour la forme quadratique
$q$. Pour un morphisme $f$ du ferm\'e compl\'ementaire, la courbe
$f(\pu)$ est contenue dans un $(\a)$-plan ou un $(\b)$-plan de $\G$,
la surface ou sa duale est alors un c\^one et son lieu singulier
abstrait n'est plus une courbe de $\ps$. Nous verrons que pour tout
\'el\'ement $f$ de $\R_d$, le lieu singulier abstrait $\P(f)$ est une
courbe de $S^2\pu$. La vari\'et\'e $\R_d$ est
munie d'action de $PGL_2$ et $PO(q)$. L'action induite de $PGL_4$
identifi\'e \`a $PSO(q)$ est l'action du groupe des automorphismes de
$\p^3$ sur les surfaces.

Nous \'etudions le morphisme $\P$ qui a un \'el\'ement de $\R_d$
associe son lieu singulier abstrait et montrons le th\'eor\`eme
suivant qui r\'epond au probl\`eme pour le lieu singulier abstrait :

\vs 0.1 cm

\noi
{\bf Th\'eor\`eme} : \textit{(\i) Pour $d\geq 3$, le morphisme $\P$ est
g\'en\'eriquement injectif modulo isomophisme et dualit\'e : la fibre
g\'en\'erale de $\P$ est une orbite sous $PO(q)$. Pour $d\geq 5$ la
fibre g\'en\'erale est isomorphe \`a $PO(q)\backsimeq
PGL_4\ltimes\{\pm 1\}$.}

\textit{(\i\i) La vari\'et\'e $\R_d$ est irr\'eductible de dimension
$4d+4$. Le morphisme $\P$ est dominant pour $d$ compris entre $3$ et
$5$. Pour $d\geq 6$, son image est rationelle, normale de dimension $4d-11$
r\'eguli\`ere en codimension $d-5$ et de degr\'e :}
\vs -0.5 cm
$$\mathrel{\mathop{\prod}\limits_{k=0}^{d-6}\frac{\binom{d+1+k}{d-5-k}
}{\binom{2k+1}{k}}}$$

\vs 0.1 cm

Pour montrer ce r\'esultat, nous utilisons l'isomorphisme exceptionnel
entre $PSO(q)$ et $PGL_4$ qui nous permet d'avoir deux points de vue
sur $\G$. Nous d\'efinissons ainsi deux morphismes naturels $\P$ et
$\Phi$ de la vari\'et\'e de $\R_d$ dans l'espace projectif des courbes
de degr\'e $d-2$ du plan. Nous montrerons que ces morphismes sont
\'egaux.

Plus pr\'ecisement, nous d\'ecrirons une stratification de la
vari\'et\'e $\R_d$ et montrerons que $\Phi$ et $\P$ co\" \i ncident
sur la plus petite strate. Nous montrerons que la vari\'et\'e
$\R_d$ est plong\'ee dans un espace projectif et que $\Phi$ et $\P$
sont alors lin\'eaires pour ce plongement. La strate
minimale engendrant cet espace, nous concluons \`a l'\'egalit\'e des
morphismes. Le morphisme $\Phi$ nous permettra de montrer
l'injectivit\'e g\'en\'erique et de v\'erifier que l'image, qui pourra
\^etre vue comme vari\'et\'e d\'eterminantielle, est normale
de lieu singulier en codimension $d-4$ et du degr\'e annonc\'e. Nous
donnerons au dernier paragraphe une param\'etrisation birationnelle de
l'image de $\P$ ce qui nous permettra d'affirmer que cette vari\'et\'e
est normale.

Par ailleurs, nous montrerons l'existence d'un ``espace de modules des
surfaces rationnelles r\'egl\'ees param\'etr\'es '' ${\cal M}(d)$
comme bon quotient
de $\R_d$ par $PGL_4$. Nous \'etudierons quelques sous-vari\'et\'es
remarquables de $\R_d$ ainsi que leurs images par $\P$. Nous
donnerons \'egalement une compactification de $\R_d$ et nous
d\'ecrirons l'image du bord. Enfin, nous \'etudierons le cas des
surfaces de degr\'e $5$ et retrouverons des r\'esultats de \cite{P1}
sur la position relative d'une cubique et d'une conique.

\vs 0.4 cm

\noi
{\bf Remerciements} : Je tiens ici \`a remercier mon directeur de
th\`ese \textsc{Laurent Gruson} pour toute l'aide qu'il m'a apport\'ee
durant la pr\'eparation de ce travail.

\vs 0.4 cm


\noi
{\bf Notations} : (\i) Notons $S_n$ la representation irr\'eductible
de dimension $n+1$ de $SL_2$. Nous identifierons $S_n$ et ${\check
S_n}$ gr\^ace \`a la forme bilin\'eaire $SL_2$-invariante sur $S_n$
(pour plus de d\'etails sur les repr\'esentation de $SL_2$, voir
\cite{SP}). Le plan $\ps$ contient une conique canonique $C_0$ qui est
l'image  du plongement de Veronese de $\p(S_1)$. Si $X$ est une courbe
de degr\'e $n$ de $\ps$, on appelera polyg\^one de Poncelet associ\'e
\`a $X$ tout polyg\^one complet dont les c\^ot\'es sont tangents \`a
$C_0$ et dont les sommets sont sur $X$. La courbe $X$ sera dite en
relation de Poncelet avec $C_0$ si elle poss\`ede au moins un
polyg\^one de Poncelet \`a $n+1$ c\^ot\'es associ\'e. On note $\crp_n$
l'ensemble de ces courbes.
On note $K_a$ le fibr\'e de Schwarzenberger conoyau de l'injection
suivante : $\S_{a-2}\oops(-1)\fl\S_a\oops$.

\noi
(\i\i) Nous noterons $G_a$ le groupe ${\rm{Aut}}(\som)/\bC^*$ des
automorphismes modulo homoth\'eties du faisceau $\som$.

\noi
(\i\i\i) Notons $V$ l'espace vectoriel $H^0\oo_{\p^3}(1)$. La forme
quadratique $q$ d\'efinissant $\G$ est donn\'ee par le produit
ext\'erieur sur $\LV$. Nous identifierons $\LV$ et $\L^2{\check V}$
gr\^ace \`a $q$. Notons $K$ (resp. $Q$) le sous-fibr\'e
(resp. quotient) tautologique de $\G$, on a la suite exacte suivante :
\vs -0.5 cm
$$0\fl K\stackrel{i}{\fl}V\ot\oo_{\g}\stackrel{\pi}{\fl}Q\fl 0$$

\begin{rem}\hskip -0.15 cm{\bf .}
{\rm (\i) Le lieu singulier abstrait $\P(S)$ d'une surface $S$ donne des
informations sur la g\'eom\'etrie du lieu singulier de $S$ : les points
pinces de la surface sont exactement les points de $\P(S)\cap
C_0$. Les points triples de $S$ ou de sa duale correspondent aux
triangles de Poncelet associ\'es \`a $\P(S)$. Plus g\'en\'eralement
les polyg\^ones de Poncelets \`a $n$ c\^ot\'es associ\'es \`a $\P(S)$
correspondent aux points $n$-uples de $S$ ou de sa duale (ceux-ci
peuvent \^etre d\'efinis gr\^ace aux techniques de \cite{ACGH} et aux
incidences points/diviseurs sur $\pu$). Enfin, la
surface $S$ est d\'eveloppable si et seulement si $\P(S)$ contient la
conique $C_0$.

(\i\i) Un \'el\'ement $f\in\Mor{d}{\g}$ d\'efinit une
surface de $\pv$ de la fa\c con suivante : la vari\'et\'e d'incidence
points/droites ${\mathbb I}$ entre $\pv$ et $\g$ est d\'efinie au
dessus de $\g$ par $\p_{\g}(Q)$ ($Q$ est le fibr\'e quotient
tautologique). Ainsi, le \textit{pull-back} par $f$ de $Q$  d\'efinit
un morphisme de la surface r\'egl\'ee $\p_{\pu}(f^*Q)$ dans ${\mathbb
I}$. Par projection, on obtient une surface rationnelle $S$ de
degr\'e $d$ de $\pv$. Si on effectue cette construction pour le
fibr\'e ${\check K}$ on obtient alors une surface rationnelle ${\check S}$ de degr\'e $d$ de
$\pvd$ qui est la surface duale de $S$. Nous dirons qu'une surface $S$
est autoduale si il existe un isomorphisme entre $\pv$ et $\pvd$ qui
identifie $S$ et ${\check S}$.

(\i\i\i) L'action de $PGL(V)$ pr\'eserve les fibr\'es tautologiques $K$
et $Q$ de la grassmannienne. Par contre, l'action de $PO(q)$ \'echange
$Q$ et ${\check K}$. Nous pouvons ainsi caract\'eriser les surfaces
autoduales comme \'etant d\'efinies par des morphismes
$f\in\Mor{d}{\g}$ tels que l'orbite de $f$ sous $PO(q)$ est la m\^eme
que celle sous $PSO(q)$. Nous verrons (proposition \ref{fibrephi}) que
ceci est \'equivalent \`a dire, pour une surface assez g\'en\'erale,
que la courbe correspondante sur $\G$ est trac\'ee sur un hyperplan de
$\plv$ non tangent \`a $\g$.}
\end{rem}

\newpage

\section{Comparaison de deux morphismes}

\subsection{Une stratification de $\R_d$}

Soit $f\in\Mor{d}{\G}$, le faisceau $f^*Q$ est localement libre de
rang $2$ et de degr\'e $d$ sur $\pu$. C'est par d\'efinition un
quotient de $V\oopu$ (gr\^ace \`a la fl\`eche $f^*\pi$) et il est
d\'ecompos\'e. Nous avons donc une identification
$f^*Q\backsimeq\opu(a)\oplus\opu(d-a)$ avec $0\leq a\leq d$. Ceci nous
am\`ene \`a donner la d\'efinition suivante :

\begin{defi}\hskip -0.15 cm{\bf .}
{\rm Un \'el\'ement $f\in\Md$ est dit de type $a$ (avec $0\leq a\leq \d$)
si $f^*Q$ est isomorphe \`a $\opu(a)\oplus\opu(d-a)$. Cette
d\'efinition est $PSO(q)$ et $PGL_2$ invariante. Notons
$\R_{d,a}$ (resp. $\R'_{d,b}$) l'ensemble des surfaces param\'etr\'ees
de type $a$ (resp. dont la duale est de type $b$).}
\end{defi}

\begin{rem}\hskip -0.15 cm{\bf .}
{\rm (\i) Le type de la duale est d\'etermin\'e par la d\'ecomposition de
${\check K}$ (et donc par celle de $K$). Cette notion n'est donc pas
invariante par dualit\'e. 

(\i\i) Une surface de $\R_{d,0}$ est un c\^one de degr\'e $d$. Les
vari\'et\'es $\R_{d,0}$ et $\R'_{d,0}$ ne sont donc pas contenues dans
$\R_d$ mais adh\'erentes \`a $\R_d$. 
Ceci nous permet d'affirmer que pour tout \'el\'ement
$f\in\R_d$, l'application lin\'eaire $V\fl H^0(f^*Q)$ est
injective (sinon la duale est de type nul).}\label{typedual}
\end{rem}

\begin{pr}\hskip -0.15 cm{\bf .}
Dans $\Md$, les vari\'et\'es $\R_{d,a}$ (resp. $\R'_{d,a}$) sont
d\'eterminantielles, de codimension $d-2a-1$ si $a<\frac{d}{2}$ et
nulle si $a=\frac{d}{2}$. Elles forment une stratification. Les deux
stratifications sont ind\'ependantes : l'intersection
$\R_{d,a}\cap\R_{d,b}'$ est de codimension $2d-2(a+b)-2$. \label{strate}
\end{pr}

\dm :
Soit $f_u$ le morphisme universel de $\pu\times\Md$ dans
$\G\times\Md$ (on note $p$ et $q$ les projections du premier produit
et $p'$ et $q'$ celles du second), la vari\'et\'e $\R_{d,a}$ est donn\'ee
par le lieu o\`u la fl\`eche de fibr\'es vectoriels :
$$R^1p_*(f_u^*{p'}^*K\otimes p^*\oo_{\pu}(-\d-1))\fl
R^1p_*(f_u^*{p'}^*V\otimes p^*\oo_{\pu}(-\d-1))$$
a un conoyau de dimension $\d-a-1$. Ceci permet de voir que ces vari\'et\'es
forment une stratification de $\Md$. De la m\^eme fa\c con on a une
stratification par le type de la duale. On peut \'egalement voir les
surfaces de type $a$ comme le lieu o\`u la fl\`eche de fibr\'es
vectoriels : 
$$R^1p_*(f_u^*{p'}^*K\otimes p^*\oo_{\pu}(-a-2))\fl
R^1p_*(f_u^*{p'}^*V\otimes p^*\oo_{\pu}(-a-2))$$
a un conoyau de dimension un. Dans ce cas la vari\'et\'e $\R_{d,a}$ a
la codimension attendue comme vari\'et\'e d\'eterminantielle (voir la
remarque \ref{dimstrate} pour la dimension de $\R_{d,a}$).

Pour d\'eterminer l'intersection de deux strates de chacune des
stratifications, on cherche les points $f\in\Md$ qui correspondent
\`a des surfaces de type $a$ dont la duale est de type $b$. Ceci
revient \`a dire que l'on a les identifications suivantes :
$f^*K\backsimeq\opu(-b)\oplus\opu(b-d)$ et $f^*Q\backsimeq\som$. Le
``pull-back'' de la suite exacte tautologique de $\g$ \`a $\pu$ nous
permet de voir $V\oopu$ comme un \'el\'ement de
$\p{\rm{Ext}}^1(\opu(a)\oplus\opu(d-a),\opu(-b)\oplus\opu(b-d))$.

R\'eciproquement, pour qu'une telle extension convienne, il faut et il
suffit que le faisceau
${\cal E}$ obtenu soit trivial. Ceci se traduit par les
\'egalit\'es $h^0{\cal E}(-1)=h^1{\cal E}(-1)=0$ ou encore par le fait
que le morphisme $\S_{a-1}\oplus\S_{d-a-1}\fl\S_{b-1}\oplus\S_{d-b-1}$
obtenu est un isomorphisme. Ceci est le cas pour un ouvert $U$ de nos
extensions. Regardons alors $U\times PGL_4$, nous avons un morphisme de ce
produit vers $\R_d$ : le faisceau ${\cal E}$ est trivial, un
\'el\'ement de $PGL_4$ nous permet d'identifier le faisceau $\som$
quotient de ${\cal E}$ \`a un quotient de $V\oopu$ ce qui nous d\'efinit un
morphisme $f$ de $\pu$ dans $\G$. Par d\'efinition, le ``pull-back''
par $f^*$ de la suite exacte tautologique de $\G$ \`a $\pu$ est
l'extention de d\'epart. La surface ainsi obtenue est de type $a$ et
sa duale est de type $b$. La fibre de ce morphisme est exactement
donn\'e par les orbites sous $G_a\times G_b$. Ceci vient du choix
d'une identification de $f^*K$ avec $\opu(-b)\oplus\opu(b-d)$ et de
$f^*Q$ avec $\som$. Il est facile de v\'erifier que l'action de
$G_a\times G_b$ est libre. La codimension de l'image, qui est
$\R_{d,a}\cap\R_{d,b}'$, est donc $2d-2(a+b)-2$.

\vs 0.3 cm

\noi
{\bf Param\'etrisation des surfaces de type $a$.}

\vs 0.1 cm

Si on a un \'el\'ement $f\in\Md$ de type $a$, alors on peut
consid\'erer la fl\`eche sur les sections $H^0(f^*\pi):V\fl
H^0(f^*Q)$ qui est, apr\`es identification de $H^0(f^*Q)$ avec
$\Som$, une application lin\'eaire de $V$ dans $\Som$ (injective si
$f$ est dans $\R_d$ d'apr\`es la remarque
\ref{typedual}). R\'eciproquement, pour se
donner une surface de type $a$, il suffit de se donner
une application lin\'eaire injective de $V$ dans $\Som$ telle que la
fl\`eche $V\ot\opu\fl\opu(a)\oplus\opu(d-a)$ soit surjective.

Il est donc naturel de consid\'erer l'ouvert $\F_a$ de
$\p{\rm{Hom}}(V,S_a\oplus S_{d-a})$ des applications lin\'eaires
injectives telles que la fl\`eche induite
$V\ot\opu\fl\opu(a)\oplus\opu(d-a)$ sur les faisceaux soit
surjective. L'ensemble des orbites de $\F_a$ sous $PGL(V)$ forme un
ouvert $\Gr_a$ de la grassmannienne $\G(4,S_a\oplus S_{d-a})$. Le
groupe $G_a$ agit sur $\F_a$ et sur $\Gr_a$. Un \'el\'ement de $\F_a$
d\'efinit un quotient localement libre de rang $2$ et de degr\'e $d$
de $V\oopu$ qui correpond \`a un morphisme de $\pu$ dans $\G$. Ceci
d\'efinit un morphisme $\Pi_a$ de $\F_a$ dans $\Md$.


\begin{rem}\hskip -0.15 cm{\bf .} 
{\rm (\i) 
L'image du morphisme $\Pi_a$ est
exactement $\R_{d,a}$. Les fibres de $\Pi_a$ sont les orbites sous
l'action de $G_a$ : il s'agit du choix de l'identification entre
$f^*Q$ et $\som$. Il est facile de v\'erifier que l'action de $G_a$
sur $\F_a$ est libre. La vari\'et\'e $\R_{d,a}$ est donc
irr\'eductible de dimension $3d+2a+5$ si $a<\frac{d}{2}$ et $4d+4$ si
$a=\frac{d}{2}$. De plus, on peut v\'erifier que l'action est propre et donc
libre, ainsi les vari\'et\'es $\R_{d,a}$ et $\R'_{d,b}$ sont lisses
(voir \cite{KO}).


(\i\i) Si on a un morphisme $u$ de $\R_d$ (ou m\^eme de
$\R_{d,a}$) vers un sch\'ema $X$ tel que $u$ est $PGL(V)$-invariant
alors on a un morphisme $u_a$ de $\Gr_a$ dans $X$ qui se factorise par
$u$.} \label{dimstrate}
\end{rem}

\subsection{L'espace des modules des surfaces}

L'injection canonique de $\G$ dans $\plv$ d\'efinit par composition
une immersion du sch\'ema $\Mor{d}{\g}$ dans $\Mor{d}{\plv}$. Ce
dernier sch\'ema est l'ouvert des \'el\'ement $\psi$ de l'espace
projectif  $\p(\hom{\LV}{S_d})$ tels que la fl\`eche $\pt:\LV\oopu\stackrel{\psi}{\fl}S_d\oopu\fl\opu(d)$ est
surjective.
Ainsi la vari\'et\'e $\R_d$ est le localement ferm\'e de
$\p(\hom{\LV}{S_d})$ d\'efini par la condition pr\'ec\'edente et le
fait que le morphisme de $\pu$ dans $\plv$ correspondant se factorise
par $\g$.

Nous \'etudions dans ce paragraphe l'action du groupe $PGL(V)$ sur
$\p(\hom{\LV}{S_d})$. Nous identifions les points stables et
semi-stables pour cette action et montrons qu'il existe un bon
quotient de $\R_d$ par $PGL(V)$. Nous l'appelerons ``espace des
modules des surfaces'' et le noterons $\MD$. Remarquons que les points
stables et semi-stables de $\p(\hom{\LV}{S_d})$ pour l'action de
$PO(q)$ sont les m\^emes que ceux pour l'action de $PGL(V)$, le
premier \'etant une extension de degr\'e $2$ du second.

\begin{pr}\hskip -0.15 cm{\bf .}
(\i) Un point $\psi$ de $\p(\hom{\LV}{S_d})$ est instable si et seulement
si il existe un vecteur $v\in V$ ou ${\check v}\in\cv$ tels que ${\rm
Im}^t\psi$ est contenu dans $v\wedge V$ ou dans ${\rm
Ker}(\LV\stackrel{\tilde{v}}{\fl}V)$ o\`u $\tilde{v}={\check
v}\wedge{\rm Id_V}-{\rm Id_V}\wedge{\check v}$.

(\i\i) Un point $\psi$ de $\p(\hom{\LV}{S_d})$ est non stable si et
seulement si il existe un vecteur $z=v\wedge w\in\LV$ tel que ${\rm
Im}^t\psi$ est contenu dans $z^{\perp}$ l'orthogonal de $z$ pour la
forme quadratique~$q$.
\label{module}
\end{pr}

\vs -0.6 cm

\dm :
Consid\'erons une base $(e_i)_{i\in [1,4]}$ de $V$ et un sous-groupe a
un param\`etre de $SL(V)$ diagonalis\'e dans cette base de telle sorte
qu'il agisse avec un poids $\a_i$ sur $e_i$ et que l'on ait
$\a_1\geq\a_2\geq\a_3\geq\a_4$ avec $\a_1+\a_2+\a_3+\a_4=0$.

Nous pouvons alors v\'erifier qu'un \'el\'ement
$x=\sum_{i<j}a_{i,j}e_i\wedge e_j$ avec $a_{i,j}\in S_d$ est instable
pour ce sous-groupe si et seulement si $a_{1,4}=a_{2,4}=a_{3,4}=0$ ou
$a_{2,3}=a_{2,4}=a_{3,4}=0$. Ceci est \'equivalent \`a dire que $x$ est
dans ${\rm Ker}(\LV\stackrel{\tilde{e_4}}{\fl}V)$ ou dans $e_1\wedge V$.

De m\^eme, l'\'el\'ement $x$ est non stable pour ce sous-groupe si et
seulement si $a_{3,4}=0$ ce qui signifie que $x$ est dans l'orthogonal
de $e_1\wedge e_2$ pour la forme quadratique $q$.

\begin{rem}\hskip -0.15 cm{\bf .}
{\rm Nous pouvons traduire g\'eom\'etriquement la proposition
pr\'ec\'edente. Soit $f$ un \'el\'ement de $\Mor{d}{\plv}$, il est instable si et seulement si
$f(\pu)$ est contenu dans un plan isotrope pour $q$. Il est non stable
si et seulement si $f(\pu)$ est contenu dans un hyperplan tangent \`a
la grassmannienne.}
\end{rem}

\begin{cor}\hskip -0.15 cm{\bf .}
Il existe un bon quotient, not\'e $\MD$, de la vari\'et\'e
$\R_d$. C'est l'espace des modules des surfaces rationnelles
r\'egl\'ees param\'etr\'ees.\label{quotmod}
\end{cor}

\dm :
Il suffit de montrer que tous les points de $\R_d$ sont
semi-stables. Mais la remarque pr\'ec\'edente nous permet de dire
qu'un point $f\in\R_d$ est semi-stable si et seulement si la courbe
$f(\pu)$ n'est pas contenue dans un $(\a)$-plan ou un $(\b)$-plan de
$\G$ ce qui signifie exactement que la surface et sa duale ne sont pas
des c\^ones.

\begin{rem}\hskip -0.15 cm{\bf .}
{\rm L'espace des modules $\MD$ est normal comme bon quotient d'une
vari\'et\'e lisse (voir \cite{MFK}). De plus sur l'ouvert des points
stables, l'action de $PGL(V)$ est libre (voir proposition
\ref{fibrephi} pour l'\'etude des orbites) et ferm\'ee (\cite{MFK}
proposition 2.4) donc (\cite{MFK} proposition 0.9) la vari\'et\'e
$\MD$ est lisse sur l'image de l'ouvert des points stables. Elle est
donc r\'eguli\`ere en codimension $2d-8$ (voir proposition
\ref{fibrephi} pour la dimension de l'image des points semi-stables).}
\end{rem}

\subsection{Le morphisme vers le lieu singulier}

Le lieu singulier d'une surface rationnelle r\'egl\'ee est l'ensemble
des points de la surface qui sont sur deux g\'en\'eratrices. Nous
d\'efinissons maintenant sch\'ematiquement le lieu singulier abstrait
d'une surface rationnelle r\'egl\'ee param\'etr\'ee. Nous utilisons
pour ceci les techniques de \cite{ACGH} et l'incidence point/diviseur
sur $\pu$. Nous v\'erifions que c'est une courbe de $S^2\pu$ sauf si
la surface ou sa duale est un c\^one.

Soit $f$ un \'el\'ement de $\Md$, alors $f^*\pi$
est une sujection de $V\oopu$ dans $f^*Q$ qui est un faisceau
localement libre de rang $2$ et de degr\'e $d$. Notons $q$ et $p$ les
projections de la vari\'et\'e d'incidence de $S^2\pu\times\pu$ sur le
premier et le second facteur, nous pouvons construire la compos\'ee
suivante : $V\oops\fl V\otimes q_*p^*\opu\fl q_*p^*(f^*Q)$. Soit $R$
le conoyau de cette fl\`eche.

\begin{defi}\hskip -0.15 cm{\bf .}
{\rm Le lieu singulier abstrait de $f\in\Md$ est le
$0^{{\rm{i\grave{e}me}}}$ id\'eal de Fitting du faisceau $R$ d\'efini
ci-dessus. Ensemblistement, il correspond aux paires de
g\'en\'eratrices qui se rencontrent, on le note $\P(f)$.}\label{deflieu}
\end{defi}

\begin{lem}\hskip -0.15 cm{\bf .}
Soit $f\in\Md$, le sch\'ema $\P(f)$ est une
courbe de degr\'e $d-2$ si et seulement si $f$ est dans
$\R_d$. L'application $\P$ ainsi d\'efinie sur $\R_d$ est
$PO(q)$-invariante et \`a valeurs dans $\p(S^{d-2}\S_2)$, la
vari\'et\'e des courbes de degr\'e $d-2$ de $\p(\S_2)$.
\end{lem}

\dm :
Le lieu singulier abstrait est une vari\'et\'e d\'eterminantielle, il
est donc vide ou de codimension au plus $1$. Il ne peut \^etre vide si
le degr\'e de la surface est au moins $3$ et si il est de codimension
nulle, ceci signifie que toutes les g\'en\'eratrices de la surface se
rencontrent. Ceci n'est possible que si la surface ou sa duale est un
c\^one.

Le faisceau $R$ admet la r\'esolution suivante :
$$H^0(f^*Q(-2))\otimes\ops(-1) \stackrel{A}{\fl}
(H^0(f^*Q)/V)\otimes\ops\fl  R \fl 0$$
Ainsi le lieu singulier abstrait est une courbe de degr\'e
$d-2$ sur $\ps$ pour tout $f\in\R_d$ et son \'equation est donn\'ee
par le d\'eterminant de $A$.

Pour montrer l'invariance sous $PO(q)$ il suffit de prouver
l'invariance par dualit\'e et par $PSO(q)\backsimeq PGL(V)$. L'invariance sous
$PGL(V)$ est \'evidente car le lieu singulier abstrait ne d\'epend que
de l'image de $V$ dans $H^0(f^*Q)$. Pour l'invariance par dualit\'e, on
dualise la suite exacte pr\'ec\'edente :
$$0\fl (H^0(f^*Q)/V {\check )}\otimes\ops(-1) \fl H^0(f^*Q(-2){\check
)}\otimes\ops \fl \e^1(R(1),\ops) \fl 0$$
En tenant compte des isomorphismes donn\'es par la dualit\'e de Serre
entre $(H^0(f^*Q)/V{\check )}$ et $H^0(f^*{\check K}(-2))$ et
entre $H^0(f^*Q(-2){\check )}$ et $H^0(f^*{\check K})/{\check V}$,
on voit que la fl\`eche d\'efinissant le lieu singulier de la duale
est la transpos\'ee de $A$. Ainsi son d\'eterminant est le
m\^eme. Cette application est donc invariante sous l'action de
$PO(q)$.

\begin{pr}\hskip -0.15 cm{\bf .}
L'application $\P$ se factorise par
$\Gr_a$ au dessus de $\R_{d,a}$, le morphisme $\P_a$ ainsi obtenu
est $PGL_2$-lin\'eaire pour le plongement de Pl\" ucker de $\Gr_a$ dans
$\p(\L^4(\Som))$.
\end{pr}

\dm :
Le morphisme $\P$ est invariant sous $PGL(V)$ il se factorise donc par
$\Gr_a$ au dessus de $\R_{d,a}$ (remarque \ref{dimstrate}). 
Sur $\Gr_a\times\ps$ le morphisme $\P_a$ est donn\'e par le d\'eterminant
de l'application suivante (on note $\V$ le sous-fibr\'e tautologique
et $\pru$ et $\prd$ les projections sur le premier et le second
facteur) :
$$(S_{a-2}\oplus S_{d-a-2})\otimes\prd^*\ops\fl
((S_a\oplus S_{d-a})/\pru^*\V)\otimes\prd^*\ops(1)$$
il est donc donn\'e par le produit ext\'erieur d'ordre $d-2$ puis 
par projection sur $\Gr_a$ par :
$$\oo_{\Gr_a}\fl (\L^4\V{\check )}\otimes S^{d-2}\S_2$$
Le morphisme $\P_a$ est donc bien lin\'eaire car $(\L^4 \V{\check )}$
est de degr\'e $1$ pour le plongement de Pl\" ucker.

\begin{rem}\hskip -0.15 cm{\bf .}
{\rm Le morphisme $\P$ se factorise par $\MD$. Le morphisme induit de
$\MD$ vers $\p(S^{d-2}\S_2)$ est compatible avec l'action de
$PGL_2$. Nous d\'ecrirons g\'eom\'etriquement ses fibres pour les
surfaces de degr\'e $5$ (voir proposition \ref{fibredeg5}).}
\end{rem}

\subsection{Le morphisme $\Phi$}


Nous construisons maintenant un second $PGL_2$-morphisme $\Phi$ de
$\R_d$ dans $\p(S^{d-2}\S_2)$ qui sera plus facile \`a
\'etudier que $\P$. Nous calculerons ses fibres et le
degr\'e de son image. Le morphisme $\Phi$ sera encore invariant sous
$PO(q)$. Nous montrerons que c'est le bon quotient de $\R_d$ sous
$PO(q)$.

Nous d\'efinissons le morphisme $\Phi$ sur tout l'espace projectif
$\p({\rm{Hom}}(\LV,\S_d))$. Rappelons que l'on a les inclusions
$\R_d\subset \Mor{d}{\p(\LV)}\subset\p({\rm{Hom}}(\LV,\S_d))$.
La seconde est donn\'ee de la fa\c con suivante : si on a un
\'el\'ement $f\in\Mor{d}{\plv}$, il se d\'ecompose en 
\begin{itemize}
\item une partie canonique : le plongement de Veronese $\pu\fl\p(S_d)$
(d\'efini par $x\mapsto x^d$) et dont l'image est $C_d$ la courbe
rationnelle normale (invariante sous $PGL_2$).
\item une projection de $\p(S_d)\dashrightarrow\plv$ qui est donn\'ee
par un \'el\'ement $\psi\in\p(\hom{\LV}{S_d})$.
\end{itemize}
L'application lin\'eaire $\psi$ est d\'efinie par les sections de la
surjection $\widetilde{\psi}:\LV\oopu\fl\opu(d)$. R\'ecipro\-quement
un \'el\'ement de
$\p({\rm{Hom}}(\LV,\S_d))$ d\'efinit un morphisme de degr\'e $d$ de
$\pu$ dans $\G$ si et seulement si l'application $\widetilde{\psi}$
ci-dessus induite sur les faisceaux est surjective.

Notons $\p(S^{2}\S_d)$ l'ensemble des formes quadratiques sur
$\S_d$. Un \'el\'ement $\psi$ de l'espace projectif
$\p({\rm{Hom}}(\LV,\S_d))$ d\'efinit 
une projection $\p(S_d)\dashrightarrow\plv$. Si on
``remonte'' la forme quadratique $q$ par cette projection, on obtient
une forme quadratique de rang inferieur ou \'egal \`a $6$ sur $\S_d$. 
Nous d\'efinissons ainsi le $PGL_2$-morphisme suivant :
\vs -0.5 cm
$$\Phi:\p{\rm{Hom}}(\LV,\S_d)\fl\p(S^2\S_d)$$
\vs - 0.9 cm
$$\psi\mapsto\psi q^{t}\!\psi$$
L'application $\Phi$ est d\'efinie en dehors du ferm\'e de
$\p{\rm{Hom}}(\LV,\S_d)$ des $\psi$ tels que l'image de$~^t\psi$ est
contenue dans un sous-espace totalement isotrope de $\LV$. Ainsi
$\Phi$ est d\'efinie sur l'ouvert des points semi-stables (voir
proposition \ref{module}) et donc sur $\R_d$ tout entier.

\begin{rem}\hskip -0.15 cm{\bf .} 
{\rm (\i) Soit $f\in\Mor{d}{\plv}$. En tant que quadrique, la forme quadratique
$\Phi(f)$ contient la courbe rationnelle normale $C_d$ si et seulement
si la courbe $f(\pu)$ \'etait trac\'ee sur $\G$. En effet, l'image de
la quadrique $\Phi(f)$ par la projection $\p(S_d)\dashrightarrow\plv$
est exactement $\G$ alors que celle de $C_d$ est
$f(\pu)$. \label{memimage}

(\i\i) L'image de $\R_d$ par $\Phi$ est donc contenue dans l'ensemble
des quadriques de $\p(S_d)$ contenant la courbe rationnelle normale
$C_d$, c'est l'espace projectif $\p(H^0\I_{C_d}(2))$. Mais l'espace
vectoriel $H^0{\cal I}_{C_d}(2)$ est le noyau de
$H^0\oo_{\p(S_d)}(2)\fl H^0\oo_{C_d}(2)$ et s'identifie \`a
$S^2\S_{d-2}$. La loi de r\'eciprocit\'e de Hermite (voir \cite {SP})
nous permet d'identifier $S^2\S_{d-2}$ \`a $S^{d-2}\S_2$. Ainsi, le
morphisme $\Phi$ d\'efini sur $\R_d$ est \`a valeur dans
$\p(S^{d-2}S_2)$ l'espace des courbes de degr\'e $d-2$ de $\ps$.}
\end{rem}

Nous allons montrer au th\'eor\`eme \ref{egalite} que ce morphisme est
le m\^eme que $\P$. L'\'egalit\'e de $\P$ et de $\Phi$ nous permettra
de d\'eterminer les fibres de $\P$ (cf. th\'eor\`eme \ref{injectif})
et le degr\'e de son image ainsi que celui de certaines
sous-vari\'et\'es particuli\`eres de cette image (cf. Th\'eor\`eme
\ref{degre}).

Notons $\Q_k$ le localement ferm\'e de
$\p(S^2\S_d)$ des formes quadratiques de rang exactement
$k$. Le ferm\'e des formes quadratiques qui contiennent la
courbe rationnelle normale $C_d$ est une sous-vari\'et\'e lin\'eaire de
$\p(S^2\S_d)$ isomorphe \`a $\p(S^{d-2}\S_2)$ (voir remarque
\ref{memimage}), notons $\q_k$ l'intersection de $\Q_k$ et de ce
sous-espace lin\'eaire. Remarquons que $\q_3$ est ferm\'e car $\aaq_2$
est vide. En effet, il n'y a pas de forme quadratique de rang au plus
$2$ contenant $C_d$ sinon cette derni\`ere serait contenue dans un
hyperplan. Notons enfin $\R^k_d$ les images r\'eciproques des $\q_k$.

\begin{lem}\hskip -0.15 cm{\bf .}
L'application $\Phi$ est d\'efinie sur les points semi-stables de
$\p\hom{\LV}{S_d}$, elle invariante sous l'action de $PO(q)$ et son
image est $\aq_6$.

La restriction de $\Phi$ \`a $\overline{\R_d}$
est d\'efinie en dehors des ferm\'es $\overline{\R_{d,0}}$ et
$\overline{\R'_{d,0}}$ et son image est $\aaq_6$. On a
$\Phi^{-1}(\Phi(\overline{\R_d}))=\overline{\R_d}$ ce qui signifie que
les fibres de $\Phi\vert_{\overline{R_d}}$ sont les m\^emes que celles
de $\Phi$ sur tout l'espace $\p\hom{\LV}{S_d}$.

L'application $\Phi$ se factorise sur $\Gr_a$ en un morphisme
lin\'eaire pour le plongement de Pl\" ucker.\label{defphi}
\end{lem}

\dm :
Nous avons d\'ej\`a vu que $\Phi$ est d\'efinie sur les points
semi-stables. Elle est \'evidement $PO(q)$-invariante. L'image de
$\p({\rm{Hom}}(\LV,\S_d))$ est contenue dans le ferm\'e des
formes quadratiques de rang inf\'erieur ou \'egal \`a
$6$. R\'eciproquement, si $w$ est une
telle forme quadratique, son noyau ${\rm Ker}w$ est de codimension au
plus $6$. Soit $N={\check S_d}/{\rm Ker}w$, c'est un espace vectoriel
de dimension au plus 6 muni d'une forme quadratique non
d\'eg\'en\'er\'ee. Si on se donne une isom\'etrie injective $i$ de $N$
dans $\L^2{\check V}$, on construit gr\^ace \`a la compos\'ee
${\check S_d}\fl N\stackrel{i}{\fl}\L^2{\check V}$ la transpos\'ee
d'un \'el\'ement de la fibre. 

Nous nous contentons pour le moment de v\'erifier que
$\overline{\R_{d,0}}$ et $\overline{\R'_{d,0}}$ sont dans le lieu base
de $\Phi$. Nous v\'erifierons au moment de la compactification de
$\R_d$ (remarque \ref{findefphi}) que le lieu base de $\Phi$ sur
$\overline{\R_d}$ est exactement donn\'e par
$\overline{\R_{d,0}}\cup\overline{\R'_{d,0}}$. L'application $\Phi$
n'est pas d\'efinie sur le lieu des points instables qui contient
 $\R_{d,0}\cup\R'_{d,0}$ (voir corollaire \ref{quotmod}) et donc son
adh\'erence. 

Un \'el\'ement de $\R_d$ est donn\'e par une
application lin\'eaire $\psi$ de $\LV$ dans $\S_d$ telle que le
morphisme $\widetilde{\psi} :\LV\oopu\fl\opu(d)$ est surjectif et que
la compos\'ee
$$\opu(-d)\stackrel{^t\!\psi}{\fl}\L^2{\check
V}\oopu\stackrel{q}{\fl}\LV\oopu\stackrel{\psi}{\fl}\opu(d)$$
est nulle. L'adh\'erence de $\R_d$ est donc donn\'ee par la seconde
condition. Par $\Phi$ on obtient alors le sous-sch\'ema de $\aq_6$ des
formes quadratiques $w$ qui v\'erifient que la compos\'ee 
$$\opu(-d)\fl\S_d\oopu\stackrel{w}{\fl}\S_d\oopu\fl\opu(d)$$
est nulle. C'est exactement $\aaq_6$.
Les fibres de la restriction de $\Phi$ \`a $\overline{\R_d}$ sont les
m\^emes que celles de $\Phi$. En effet, la condition ``la courbe
est dans $\G$'' se traduit exactement par le fait que la quadrique
image contient $C_d$ ce qui se voit sur l'image.

Il nous reste \`a montrer que l'application induite sur $\Gr_a$ (par
invariance sous $PGL(V)$) est lin\'eaire pour le plongement de
Pl\" ucker. Si $\V$ est le sous-fibr\'e tautologique et $\psi_u$ est
le morphisme universel de $\L^2\V$ dans $S_d\ot\oo_{\Gr_a}$, le
morphisme $\Phi$ est d\'efini de la fa\c con suivante : 
$$\L^4\V\stackrel{q}{\fl}S^2(\L^2\V)\stackrel{s^2\psi_u}{\fl} S^2\S_d\otimes\oo_{\Gr_a}$$
qui est bien lin\'eaire pour le plongement de Pl\" ucker car $\L^4\V$
est de degr\'e $-1$.

\begin{pr}\hskip -0.15 cm{\bf .}
Les fibres de $\Phi$ sont isomorphes \`a $PO(q)$ au dessus de
$\Q_6$ et  $PSO(q)$ au dessus de $\Q_5$. Dans ces deux cas la fibre est une
orbite ferm\'ee sous $PO(q)$.
Au dessus de $\Q_4$, les fibres de $\Phi$ ont deux composantes
irr\'eductibles de dimension $d+11$ qui se coupent en un ferm\'e de
dimension $14$ qui est l'unique orbite ferm\'ee de cette fibre. Au
dessus de $\Q_3$, elles sont irr\'eductibles de dimension $d+11$ et
contiennent une unique orbite ferm\'ee de dimension $12$.

L'image de $\R_d$ est exactement $\aaq_6$.\label{fibrephi}
\end{pr}

\dm :
Si $G$ est un sous-groupe d'un groupe lin\'eaire, on note $PG$ le
quotient du groupe $G$ par les homoth\'eties.

Soit $w$ une forme quadratique de $\S_d$. Si $\psi\in \Phi^{-1}(w)$, alors
${\rm{Ker}}^t\!\psi$ est inclus dans ${\rm{Ker}}w$. Notons
$I(\psi)=^t\!\!\!\psi({\rm{Ker}}w)$ qui est un sous-espace vectoriel
isotrope de $\LV$ et $N(w)$ l'espace $S_d/{\rm {Ker}}(w)$ (qui est
muni d'une forme quadratique non d\'eg\'en\'er\'ee induite par
$w$). La suite exacte :
\vs -0.4 cm
$$0\fl I(\psi)\fl {\rm{Im}}^t\!\psi\fl N(w)\fl 0$$
\vs -0.1 cm
impose que $I(\psi)$ est exactement le noyau de la restriction de $q$
\`a ${\rm{Im}}(^t\!\psi)$. Ainsi, 
${\rm{Im}}(^t\!\psi)$ est n\'ecessairement contenu dans $I(\psi)^{\perp}$. Les rangs
$r(w)$ et $r(\psi)$ de $w$ et $\psi$ et la dimension $i(\psi)$ de
$I(\psi)$ v\'erifient les deux conditions suivantes : 
\vs -0.4 cm
$$r(\psi)=r(w)+i(\psi)\ \ {\rm{et}}\ \ r(\psi)\leq 6-i(\psi)$$
\vs -0.1 cm
Une fois fix\'es les entiers $r(\psi)$ et $i(\psi)$ v\'erifiant ces
conditions, la fibre au dessus de $w$ est donn\'ee par les choix
suivants : un sous-espace vectoriel ${\rm{Ker}}^t\!\psi$ de
${\rm{Ker}}w$ de codimension $i(\psi)$ ; un sous-espace
${\rm{Im}}(^t\!\psi)$ de dimension $r(\psi)$ de $\LV{\check{~}}$ tel
que la restriction de la forme quadratique $q$ a pour noyau $I(\psi)$
un espace de dimension $i(\psi)$ ; une isom\'etrie \`a homoth\'etie
pr\`es de $\S_d/{\rm{Ker}}(^t\!\psi)$ dans
${\rm{Im}}(^t\!\psi)$. Le groupe $PO(q)$ agit transitivement sur les
couples form\'es d'objets des deux derniers types.

Les valeurs de $r(w)$, $r(\psi)$ et $i(\psi)$ possibles sont alors (on
se contente des cas $r(w)\geq 3$) : 
\begin{itemize}
\item si $r(w)=6$, alors $r(\psi)=6$ et
$i(\psi)=0$ 
\item si $r(w)=5$, alors $r(\psi)=5$ et $i(\psi)=0$
\item si $r(w)=4$, alors $r(\psi)=5$ et $i(\psi)=1$ ou $r(\psi)=4$ et
$i(\psi)=0$
\item si $r(w)=3$, alors $r(\psi)=4$ et $i(\psi)=1$ ou $r(\psi)=3$ et
$i(\psi)=0$. 
\end{itemize}

On commence par traiter les cas o\`u $i(\psi)$ est nul. Dans ce cas, la
fibre est donn\'ee par un sous-espace ${\rm{Im}}(^t\!\psi)$ de
$\LV{\check{~}}$ de dimension $r(w)$ tel que la restriction de
la forme quadratique $q$ est non d\'eg\'en\'er\'ee et une isom\'etrie
de $\S_d/{\rm{Ker}}(w)$ dans ${\rm{Im}}(^t\!\psi)$. Le th\'eor\`eme
de Witt nous permet de dire que tous les sous-espaces
${\rm{Im}}(^t\!\psi)$ sont conjugu\'es sous $PO(q)$, la fibre est
ainsi une orbite sous ce groupe. Elle est toujours ferm\'ee. Le
stabilisateur est alors $PO({\rm{Im}}(^t\!\psi)^{\perp})$. Ceci
d\'ecrit les cas de rang $5$ et $6$ et les ferm\'es de dimensions $14$
et $12$ des cas de rang $4$ et $3$. Dans le cas o\`u la forme
quadratique est de rang $5$, la fibre est isomorphe \`a $PSO(q)$ est
c'est une orbite sous $PO(q)$. Les actions de $PO(q)$ et de $PSO(q)$
pour les points de $\Phi^{-1}(\Q_5)$ sont les m\^emes ce qui signifie
que les surfaces de $\R_d^5$ sont autoduales. Ces courbes sont exactement
les courbes trac\'ees sur un hyperplan r\'egulier (non tangent \`a la
quadrique) de $\p(\LV)$. De m\^eme sur les deux ferm\'es de
dimensions 14 et 12, l'action de $PSO(q)$ est la m\^eme que celle de
$PO(q)$ ainsi ces surfaces sont autoduales.

Si $i(\psi)$ vaut $1$ (dans les cas $r(w)=4$ ou $r(w)=3$), la fibre est
alors donn\'ee par le choix d'un hyperplan de ${\rm{Ker}}w$ (c'est
$\p({\rm{Ker}}w)$ qui est de dimension $d-4$ ou $d-3$), le choix d'un
sous-espace $W$ de dimension $5$ ou $4$ de $\LV{\check{~}}$ tel que la
restriction de la forme quadratique $q$ \`a $W$ a un noyau de
dimension $1$ (ce choix correspond \`a $PO(q)/{\rm{Stab}}(W)$ qui est de dimension $4$ ou
$7$) et le choix d'une isom\'etrie \`a homoth\'etie pr\`es (c'est
$PO(W)$ qui est isomorphe \`a $PO(4)\times {\bC}^*\times{\bC}^4$
ou $PO(3)\times {\bC}^*\times{\bC}^3$ et est de dimension $11$ ou
$7$). On peut remarquer qu'une fois le choix de l'\'el\'ement de
$\p({\rm{Ker}}w)$ r\'ealis\'e, la fibre est une orbite sous $PO(q)$. En
effet, pour $r(w)=4$, on a $PO(W)={\rm{Stab}}(W)$ donc la fibre est
isomorphe \`a $\p({\rm{Ker}}w)\times PO(q)$ et a deux composantes
connexes de dimension $d+11$. Dans le cas $r(w)=3$, on a l'\'egalit\'e
${\rm{Stab}}(W)=PO(W)\times^{PO(W)\cap PO(W^{\perp})}PO(W^{\perp})$ o\`u
$PO(W^{\perp})={\bC}^*\times{\bC}\times\{{\rm{Id}},\sigma_H\}$
o\`u $H$ est un hyperplan contenant $W$ et $PO(W)\cap
PO(W^{\perp})={\bC}^*$. On voit alors que la fibre est isomorphe \`a
$\p({\rm{Ker}}w)\times (PO(q)/{\bC}\times\{{\rm{Id}},\sigma_H\})$. Le
sous-groupe ${\bC}\times\{{\rm{Id}},\sigma_H\}$ de $PO(W^{\perp})$
est form\'e des \'el\'ements dont la restriction \`a $I(\psi)$ est
l'identit\'e. On voit ainsi que la fibre est isomorphe \`a
$\p({\rm{Ker}}w)\times PSO(q)/{\bC}$ et qu'elle est irr\'eductible
de dimension $d+11$.

Il reste \`a v\'erifier que les parties correspondant \`a $i(\psi)$ nul
sont adh\'erentes \`a celles o\`u $i(\psi)=1$. Pour le voir, on choisit
$w$ de rang $4$ ou $3$ et $\psi\in \Phi^{-1}(w)$ tel que $i(\psi)=0$ et on
construit une d\'eformation $\psi_t$ d'\'el\'ement tels que $i(\psi_t)$
est non nul qui tend vers $\psi$ quand $t$ tend vers
$0$. L'application $\psi$ est donn\'ee par une injection $\psi
:N(w)\fl\LV$. Soit $K$ un hyperplan de ${\rm{Ker}}(w)$ et
$N=\S_d/K$. Soit $I$ le noyau de dimension $1$ de $N\fl N(w)$, $x$ un
g\'en\'erateur de $I$ et $y$ un \'el\'ement isotrope de
$(\psi(N(w)))^{\perp}$, on d\'efinit l'application $\psi_t$ sur
$N=I\oplus N(w)$ par $\psi_t\vert_{N(w)}=\psi$ et $\psi_t(x)=ty$. Elle
convient.

Dans le cas $r(w)=3$ et $i(\psi)=0$, le noyau
${\rm{Ker}}(^t\!\psi)={\rm{Ker}}w$ rencontre la courbe rationnelle
normale $C_d$ en $d-2x$ points et se projette sur une conique avec
multiplicit\'e $x$. Ainsi, la courbe image est de degr\'e $2x$. On est
plus dans
$\R_d$ sauf si $d$ est pair et $x=\d$. On verra au paragraphe
\ref{bord} que les autres cas
sont adh\'erents \`a $\R_d$. Si par contre $i(\psi)=1$, alors on peut
trouver un \'el\'ement $\psi$ qui nous donne une courbe de degr\'e
$d$. De m\^eme, au dessus des formes quadratiques de rang $4$, $5$ et
$6$, on peut trouver des \'el\'ements de $\R_d$. Pour ceci, il suffit
de choisir un sous-espace ${\rm{Ker}}^t\!\psi$ de ${\rm{Ker}}w$ ne
rencontrant pas $C_d$. On voit ainsi que $\q_i$ est l'image de
$\R_d^i$.

\begin{rem}\hskip -0.15 cm{\bf .}
{\rm (\i) Les surfaces de $\R^5_d$ correspondent aux courbes trac\'ees
sur un hyperplan non tangent \` a $\G$. Les surfaces de $\R^4_d$
correspondent aux courbes trac\'ees sur un hyperplan tangent \`a $\G$
(cas g\'en\'eral) ou \`a celles trac\'ees sur un espace projectif de
dimension $3$ de $\p(\LV)$ (ferm\'e de codimension $d-3$). 
Enfin, les surfaces de $\R_d^3$ correspondent aux courbes trac\'ees
sur un espace projectif de dimension $3$ tangent \`a $\G$.

(\i\i) Les fibres de l'application de ${\cal M}(d)$ vers
$\p(S^2\S_{d})$ d\'eduite de $\Phi$ sont form\'ees par deux points au
dessus de $\q_6$, un seul point au dessus de $\overline{\q_5}$.}
\end{rem}

\newpage

\begin{pr}\hskip -0.15 cm{\bf .}
Le morphisme $\Phi$ est le bon quotient de l'ouvert des points
semi-stables de $\p(\hom{\LV}{S_d})$ par $PO(q)$. En particulier, le
morphisme $\Phi$ resteint \`a $\R_d$ est un bon quotient pour l'action
de $PO(q)$.
\end{pr}

\dm :
Les points stables et semi-stables de $\p(\hom{\LV}{S_d})$ pour
l'action de $PO(q)$ sont les m\^emes que ceux pour l'action de
$PGL(V)$. Ils sont d\'ecrits \`a la proposition \ref{module}. Ainsi
l'ouvert de d\'efinition de $\Phi$ est exactement l'ouvert des
points semi-stables. Par ailleurs $\Phi$ est $PO(q)$ invariant donc il
se factorise par le bon quotient sur cet ouvert. Enfin, l'image de
$\Phi$ est $\overline{\Q}_6$ qui est de cohen-macaulay 
et lisse en dehors de $\overline{\Q}_5$ (voir \cite{ACGH}). Elle est
donc r\'eguli\`ere en codimension $d-5$ et donc en
codimension $1$ pour $d\geq 6$ (pour $d\leq 5$ elle est lisse). Le
crit\`ere de Serre nous permet de dire que $\overline{\Q}_6$ est normale. Par ailleurs, la proposition
\ref{fibrephi} nous dit qu'il y a une unique orbite ferm\'ee dans
chaque fibre de $\Phi$. Ainsi le morphisme du quotient
$\p(\hom{\LV}{S_d})^{ss}/PO(q)$ vers $\overline{\Q}_6$ est une
bijection ensembliste. La normalit\'e de $\overline{\Q}_6$ nous permet
d'appliquer le th\'eor\`eme principal de Zariski pour affirmer que
c'est un isomorphisme.

\begin{cor}\hskip -0.15 cm{\bf .}
La vari\'et\'e $\aaq_6$ est normale r\'eguli\`ere en codimension
$d-5$.
\end{cor}

\dm :
Nous savons que $\aaq_6$ est l'image de $\Phi$ par $\R_d$. C'est donc
le bon quotient d'une vari\'et\'e lisse d'o\`u la normalit\'e. Par
ailleurs, sur l'ouvert $\R_d^6$, l'action est propre (les points sont
stables) et libre (proposition \ref{fibrephi}). Ainsi la vari\'et\'e
$\q_6$ est lisse et son compl\'ementaire $\aaq_5$ qui est donc exactement le lieu
singulier de $\aaq_6$ est en codimension $d-4$. 

\subsection{Comparaison des deux morphismes}

Nous montrons dans ce paragraphe le th\'eor\`eme suivant :

\begin{theo}\hskip -0.15 cm{\bf .}
Les morphismes $\P$ et $\Phi$ sont
\'egaux.\label{egalite}
\end{theo}

Le principe de la d\'emonstration est d'\'etudier les applications
lin\'eaires $\Phi_a$ et $\P_a$ de $\Gr_a$ dans
$\p(S^{d-2}\S_2)$. Elles se prolongent sur le plongement de Pl\" ucker
en deux applications $PGL_2$-lin\'eaires que nous notons encore
$\Phi_a$ et $\P_a$ de 
$\Lambda^4(\S_a\oplus\S_{d-a})$ dans $S^{d-2}\S_2$. Il suffit de
v\'erifier qu'elles correspondent au m\^eme quotient. 

Nous \'etudions
ces morphismes sur la strate minimale (celle des surfaces de type $1$)
et nous montrons qu'ils co\"\i ncident : les deux morphismes se
factorisent par $\L^2\S_{d-1}\ot\L^2\S_1$.

Nous montrons ensuite que si $\P$ et $\Phi$ coincident sur cette strate alors
ils sont \'egaux. Par par dualit\'e, les morphismes $\Phi$ et $\P$
sont \'egaux sur les surfaces dont la duale est de type $1$. Par
lin\'earit\'e il suffit de montrer que la strate des surfaces dont la
duale est de type $1$ engendre tout l'espace vectoriel.

\begin{lem}\hskip -0.15 cm{\bf .}
Les morphismes $\P_1$ et $\Phi_1$ sont \'egaux, ils correspondent \`a la
projection canonique de $\L^4(\S_1\oplus\S_{d-1})$ vers
$\L^2\S_{d-1}\ot\L^2\S_1$.
\end{lem}

\dm :
Nous commen\c cons par \'etudier $\P_a$ : consid\'erons le morphisme :
$$(\S_{a-2}\oplus\S_{d-a-2})\oops\stackrel{m}{\fl}(\S_a\oplus\S_{d-a})
\oops(2)$$
produit des deux multiplications. Le morphisme $\P$ est donn\'e par
$H^0(\L^{d-2}m)$. Il est donc d\'efini par le produit tensoriel des
deux applications $SL_2$-lin\'eaires
$\L^{a-1}S_{a-2}\fl\L^{a-1}S_a\ot S^{a-1}S_2$ et
$\L^{d-a-1}S_{d-a-2}\fl\L^{d-a-1}S_{d-a}\ot S^{d-a-1}S_2$
Nous pouvons donc \'ecrire $\P_a$ comme \'etant l'application
lin\'eaire suivante (en tenant compte des isomorphismes entre
$\L^{k+1}S_n$ et $\L^{n-k}{\check S_n}$) :
$$\L^4(\Som)\fl\L^{2}S_a\ot\L^{2}S_{d-a}\fl S^{d-2}S_2$$
	

\'Etudions maintenant $\Phi_a$ : si on a une application $\varphi$ de $V$ vers 
$\S_a\oplus\S_{d-a}$ alors la fl\`eche de $\LV$ vers $S_d$ est donn\'ee par 
l'application $\pi :
\L^2(\S_a\oplus\S_{d-a})\fl\S_a\ot\S_{d-a}\fl\S_d$ compos\'ee avec
$\L^2\varphi$. Ainsi, le morphisme$\Phi_a$ est donn\'e par 
$$\L^4(\S_a\oplus\S_{d-a})\stackrel{^t\alpha}{\fl}
S^2(\L^2(\S_a\oplus\S_{d-a})) \fl S^2(S_a\ot S_{d-a})\fl S^2{S}_d$$
Mais le morphisme $S^2(S_a\ot S_{d-a})\fl \L^4(\S_a\oplus\S_{d-a})$ se
factorise par $\L^{2}S_a\ot\L^{2}S_{d-a}$ :
$$\Big(\sum_i(x_i\ot y_i)\Big).\Big(\sum_j(x'_j\ot
y'_j)\Big)\mapsto
-\sum_{i,j}(x_i\wedge x'_j,0)\wedge(0,y_i\wedge y'_j)$$
\vs -0.3 cm
Cette fl\`eche de $\L^4(\S_a\oplus\S_{d-a})$ dans
$\L^{2}S_a\ot\L^{2}S_{d-a}$ est la m\^eme que celle de $\P_a$. Il nous
suffit de v\'erifier que les deux applications $SL_2$-lin\'eaires de
$\L^{2}S_a\ot\L^{2}S_{d-a}$ dans $S^{d-2}S_2$ sont les m\^emes. 

Mais dans le cas $a=1$, le terme $\L^2\S_{d-1}\ot\Lambda^2\S_1$ est
isomorphe \`a $S^{d-2}\S_2$ et \`a $S^2S_{d-2}$. Ainsi dans ce cas les
deux morphismes correspondent au quotient
$\L^4(\S_1\oplus\S_{d-1})\fl\L^{2}S_1\ot\L^{2}S_{d-1}$ et apr\`es
identification de $S^{d-2}S_2$ et $S^2S_{d-2}$ \`a $\L^{2}S_{d-1}$ on
a l'\'egalit\'e entre $\P_1$ et $\Phi_1$.

\vs 0.4 cm

Les morphismes $\P$ et $\Phi$ \'etant invariants par dualit\'e nous
savons qu'ils co\"\i ncident sur le ferm\'e de $\Gr_{\d}$ des surfaces
de type g\'en\'eral dont la duale est de type $1$. Mais $\P_{\d}$ et
$\Phi_{\d}$ sont lin\'eaires pour le plongement de Pl\" ucker de
$\Gr_{\d}$, il suffit donc de montrer que ce ferm\'e n'est pas contenu
dans un hyperplan et est r\'eduit. La proposition suivante nous permet
de conclure :

\begin{pr}\hskip -0.15 cm{\bf .}
Les surfaces de type g\'en\'eral dont la duale est de type $1$ ne sont
pas sur un hyperplan du plongement de Pl\" ucker de $\Gr_{\d}$ et
forment un sch\'ema r\'eduit.
\end{pr}

\dm :
Nous savons que cette sous-vari\'et\'e est d\'eterminantielle
(proposition \ref{strate}). Nous pouvons donc calculer la r\'esolution
de son id\'eal gr\^ace au complexe d'Eagon-Northcott 
associ\'e. En effet, la vari\'et\'e est donn\'ee par la non
surjectivit\'e de la fl\`eche suivante sur $\Gr_{\d}$ (nous avons
not\'e $\V$ le sous-fibr\'e tautologique de $\Gr_{\d}$) :
$$\V\ot\S_{d-3}\fl(\S_{\d+d-3}\oplus\S_{2d-\d-3})\ot\oo_{\Gr_{\d}}$$
Notons $L$ (resp. $M$) le faisceau de gauche (resp. droite) et $l$
(resp. $m$) son rang. La r\'esolution son id\'eal ${\cal I}$ est
donn\'ee par (voir par exemple \cite{GP}) : 
$$0\fl\L^m {\check M}\ot\L^{l}L\ot S^{l-m}({\check M})\fl\cdots\fl
\L^{m}{\check M}\ot\L^{m}L\fl {\I}\fl 0$$
Il nous reste \`a calculer les sections globales de ${\I}(1)$ donc la 
cohomologie des termes de ce complexe. Pour voir que $H^0\I(1)$ est nul il 
nous suffit de montrer que les groupes de cohomologie 
$H^{l-m-p}(\L^m {\check M}\ot\L^{l-p}L\ot S^{l-m-p}({\check M})(1))$
sont nuls. Or seul le terme $L$ est non trivial, il suffit donc de
montrer que les groupes $H^{l-m-p}(\L^{l-p}L(1))$ sont nuls. Nous
pouvons calculer ce faisceau gr\^ace aux foncteurs de Schur (\cite{FH}) :
\vs -0.5 cm
$$\L^{l-p}L=\bigoplus_{\lambda}(\S_{\lambda}\V\ot\S_{\lambda'}\S_{d-3})$$
il reste donc \`a calculer la cohomologie de $\S_{\lambda}\V$ pour les 
partitions qui apparaissent et on trouve le r\'esultat d'annulation (la 
cohomologie des foncteurs de Schur des fibr\'es tautologiques est 
connue, voir \cite{DE}). Ici les $\S_{\lambda}\V$ sont \`a cohomologie
compl\`etement nulle.

Nous devons maintenant \'eliminer le cas o\`u les courbes de type $1$
seraient sur un hyperplan \'epaissi. Ce n'est pas le cas : les
surfaces de $\Gr_{\d}$ dont la duale est de type $1$ est l'image
de l'incidence entre $\Gr_1$ et $\Gr_{\d}$ donn\'ee par la
dualit\'e. Elle est donc r\'eduite.

\begin{theo}\hskip -0.15 cm{\bf .}
Le morphisme $\P$ est g\'en\'eriquement injectif modulo isomorphisme et
dualit\'e.\label{injectif}
\end{theo}

\dm :
Il suffit de combiner le th\'eor\`eme \ref{egalite} et la proposition
\ref{fibrephi}.

\section{Applications}

Nous allons maintenant utiliser le r\'esultat du th\'eor\`eme
\ref{egalite} pour \'etudier l'image de $\P$ et celle de certaines
sous-vari\'et\'es remarquables de $\R_d$. Nous \'etudierons
\'egalement la compatification de $\R_d$ dans $\p\hom{\LV}{S_d}$ ainsi
que l'image du bord par $\P$. Nous donnerons enfin quelques
applications g\'eom\'etriques, en particulier l'exemple des surfaces
de degr\'e $5$ nous permettra de d\'ecrire les positions d'une cubique
plane par rapport \`a une conique.

\subsection{\'Etude des $\R_d^k$}

Nous montrons dans ce paragraphe le th\'eor\`eme suivant :

\begin{theo}\hskip -0.15 cm{\bf .}
Les vari\'et\'es $\R_d$ et $\R^5_d$ sont irr\'eductibles et lisses
de dimensions $4d+4$ et $3d+8$.
La vari\'et\'e $\R_d^4$ a $d-1$ composantes irr\'eductibles de
dimension $3d+7$.
La vari\'et\'e $\R_d^3$ a $\d$ composantes irr\'eductibles de
dimension $2d+9$.
\label{compirr}
\end{theo}

\begin{pr}\hskip -0.15 cm{\bf .}
Les vari\'et\'es $\R_d$ et $\R^5_d$ sont irr\'eductibles et lisses
de dimensions $4d+4$ et $3d+8$.\label{irredrd}
\end{pr}

\dm :
La vari\'et\'e $\R_d$ est un ouvert du sch\'ema des morphismes de
degr\'e $d$ de $\pu$ dans la grassmannienne $\G$. Les r\'esultats de
\cite{P2} nous permettent de conclure \`a la lissit\'e,
l'irr\'eductibilit\'e et la dimension de cette vari\'et\'e. 

Nous avons un morphisme propre de la vari\'et\'e $\R^5_d$ vers les
hyperplans de $\p(\LV)$ non tangents \`a $\G$ : \`a $f\in\R_d^5$ on associe
l'hyperplan sur lequel l'image est trac\'ee (il est unique sinon
par $\Phi$ on aurait une quadrique de rang au plus $4$). La fibre de
ce morphisme est le sch\'ema $\Mor{d}{Q_3}$ des morphismes de degr\'e
$d$ de $\pu$ dans $Q_3$ ``la'' quadrique non singuli\`ere de $\p^4$
(cette quadrique est la quadrique de rang $5$ d\'ecoup\'ee dans $\G$
par l'hyperplan, elle est isomorphe \`a $Q_3$). Les r\'esultats de
\cite{P2} nous permettent d'affirmer que la fibre est irr\'eductible,
lisse de dimension constante \'egale \`a $3d+3$ d'o\`u le
r\'esultat.

\vs 0.2 cm

\def \U {{\bf U}}

\noi
{\bf La vari\'et\'e $\R_d^4$}

\vs  0.1 cm

Nous allons maintenant d\'ecrire les composantes irr\'eductibles de
$\R_d^4$. Nous avons vu \`a la proposition \ref{fibrephi} qu'un
ouvert dense $\U_4$ de $\R_d^4$ est form\'e des \'el\'ements $f$
tels que la courbe $f(\pu)$ est contenue dans un unique
hyperplan tangent \`a $\G$ (cas $i(f)=1$). Il nous suffit donc de
d\'eterminer les composantes irr\'eductibles de l'ouvert $\U_4$.

Nous avons alors un morphisme propre $\Upsilon_4:\U_4\fl\G$ qui a $f$ associe
le point de contact de l'unique hyperplan tangent \`a $\G$ contenant
$f(\pu)$. Il s'agit maintenant d'\'etudier la fibre de
$\Upsilon_4$. Cette fibre au dessus de $L_0\in\G$ est isomorphe \`a
$\Mor{d}{\cC_{L_0}}$ o\`u $\cC_{L_0}$ est le c\^one de $\G$ form\'e
par les droites qui rencontrent $L_0$.

\def \cclo {\cC_{L_0}}
\def \qlo {Q_{L_0}}
\def \ecl {\widetilde{\cC}_{L_0}}
\def \N {\mathbb{N}}

Notons $\widetilde{\cC}_{L_0}$ la vari\'et\'e : $\{(P,L,H)\in\pv\times
\G\times\pvd\ /\ P\in L_0\subset H\ \ \ {\rm{et}}\ \ \ P\in L\subset
H\}$, c'est l'\'eclatemement $\pi$ de $\cC_{L_0}$ au sommet du
c\^one. Nous avons une projection $p$ de $\ecl$ vers la vari\'et\'e $\qlo$
suivante : $\{(P,H)\in\pv\times\pvd\ /\ P\in L_0\subset H\}$. Cette
vari\'et\'e est une quadrique isomorphe \`a $L_0\times{\check
L_0}$. 

Il est facile de v\'erifier que $\ecl$ est la fibration en droites
projectives au dessus de $\qlo$ associ\'ee au faisceau
$\oo_{L_0\times{\check L_0}}(-1,0)\oplus\oo_{L_0\times{\check
L_0}}(0,1)$. Le groupe de Picard de $\ecl$ est donc $\Z^3$ (nous
prendrons pour base des cycles les degr\'es sur $\qlo$ et le degr\'e
relatif pour $p$).

\begin{lem}\hskip -0.15 cm{\bf .}
Pour tout $\gamma=(a,b,c)\in\N^3$, le sch\'ema $\Mor{\gamma}{\ecl}$ est
irr\'eductible et lisse de dimension $2a+2b+c+3$.
\end{lem}

\dm :
Le th\'eor\`eme principal de \cite{P2} nous permet de dire que le
sch\'ema des morphismes $\Mor{(a,b)}{\qlo}$ est irr\'eductible et
lisse de
dimension $2(a+b+1)$. Par ailleurs, la proposition 4 de \cite{P2}
nous permet de dire que comme le degr\'e relatif est positif, alors
$\Mor{\gamma}{\ecl}$ est irr\'eductible et lisse de dimension
$2(a+b+1)+c+1$.

\vs 0.4 cm

L'\'eclatement $\pi$ d\'efinit un morphisme
$\pi:\Mor{\gamma}{\ecl}\fl\Mor{\pi_*\gamma}{\cclo}$. De plus si
$f$ est un \'el\'ement de $\Mor{(a,b,c)}{\ecl}$ alors $f(\pu)$
rencontre le diviseur exceptionnel en $\frac{1}{2}(c-(a+b))$ points,
le degr\'e de son image dans $\cclo$ est donc
$\frac{1}{2}(a+b+c)$. Ainsi $\pi(f)$ est de degr\'e $d$ si et
seulement si $c=2d-(a+b)$.

L'image de $\Mor{(a,b,2d-(a+b))}{\ecl}$ dans $\Mor{d}{\cclo}$ est
contenu dans l'ensemble des morphismes dont l'image passe
$d-(a+b)$ fois par le sommet du c\^one. Si $d<a+b$ ceci
impose que l'image du morphisme est contenu dans le diviseur
exceptionnel, il se contracte donc sur le sommet par $\pi$. Nous
supposerons desormais que $d\geq a+b$. R\'eciproquement, si on a un
morphisme de $\pu$ dans $\cclo$ qui passe $d-(a+b)$ fois par le
diviseur exceptionnel, sa transform\'ee stricte dans
$\ecl$ est dans le sch\'ema $\Mor{(a,b,2d-(a+b))}{\ecl}$ (ou
$\Mor{(b,a,2d-(a+b))}{\ecl}$) et on
retrouve le morphisme de d\'epart par projection par $\pi$. Nous avons
donc des immersions
$\pi:\Mor{(a,b,2d-(a+b))}{\ecl}\fl\Mor{d}{\cclo}$.

\begin{pr}\hskip -0.15 cm{\bf .}
Pour $a$ fix\'e, les vari\'et\'es $\pi(\Mor{(a,b,2d-(a+b))}{\ecl})$
avec $d>a+b$ sont dans l'adh\'erence de
$\pi(\Mor{(a,d-a,d)}{\ecl})$.\label{adhe}
\end{pr}

\dm :
Notons $W$ l'espace vectoriel $H^0\oo_{L_0}(1)$, c'est un quotient de
rang $2$ de $V$. Notons $N_W$ le noyau de la surjection de $V$ dans
$W$. Si $f(\pu)$ est dans le c\^one $\cclo$, alors l'image de la
fl\`eche $N_W\oopu\fl f^*Q$ est de rang 1. De plus le
conoyau de cette  fl\`eche est de rang $2$ aux points de $\pu$ qui
s'envoient sur le sommet du c\^one et de rang $1$ ailleurs. 

Consid\'erons donc un morphisme $f\in\Upsilon_4^{-1}(L_0)$ tel que le
morphisme $N_W\oopu\fl f^*Q$ ait pour conoyau
$\opu(a)\oplus\oo_D$ o\`u $D$ est un diviseur de degr\'e
$x$. L'image de $N_W\oopu$ dans $f^*Q$ est donc
$\opu(d-(a+x))$ et le faisceau $\opu(a)\oplus\oo_D$ est le conoyau de
$\opu(-(a+x))\stackrel{r_0}{\fl} W\oopu$.

En consid\'erant les fl\`eches duales (correspondant \`a la surface
duale), nous avons un morphisme de $\cv\oopu$ dans $f^*{\check K}$. La
fl\`eche $^tr_0$ de ${\check W}\oopu$ dans $\opu(a+x)$ a pour image
$\opu(a)$ et conoyau $\oo_D$. Le morphisme ${\check W}\oopu\fl
f^*{\check K}$ se factorise par $\opu(a+x)$ et a pour conoyau
le faisceau $\opu(d-(a+x))\oplus\oo_D$. Ainsi, il suffit de se donner
une famille $r_{\lambda}$ de fl\`eches de ${\check W}\oopu$ dans $\opu(a+x)$
sujective au point g\'en\'erique et \'egale \`a $^tr_0$ au point
sp\'ecial pour constuire une famille de morphismes dont l'image ne
passe pas par le sommet au point g\'en\'erique et qui vaut $f$ au
point sp\'ecial.

\begin{cor}\hskip -0.15 cm{\bf .}
Les composantes irr\'eductibles de $\Mor{d}{\cclo}$ sont les
adh\'erences des images $\pi(\Mor{(a,d-a,d)}{\ecl})$ pour $0\leq a\leq
d$. Elles sont toutes de dimension $3d+3$.
\end{cor}

Notons $\R_d^4(a)$ le ferm\'e de $\U_4$ des morphismes $f$ qui sont
dans la composante de $\Upsilon_4^{-1}(\Upsilon_4(f))$ contenant
$\pi(\Mor{(a,d-a,d)}{\widetilde{\cC}_{\Upsilon_4(f)}})$.

\begin{cor}\hskip -0.15 cm{\bf .}
Les vari\'et\'es $\R_d^4(a)$ pour $1\leq a\leq d-1$ sont les
composantes irr\'eductibles de $\R_d^4$.  Elles sont de dimension
$3d+7$. Les composantes $\R_d^4(a)$ et $\R_d^4(d-a)$ sont  duales
l'une de l'autre.\label{compr4}
\end{cor}

\dm :
Les vari\'et\'es $\R_d^4(0)$ et $\R_d^4(d)$ ne sont pas dans $\R_d^4$
(mais dans son adh\'erence) car dans ce cas la courbe est trac\'ee sur
un $(\beta)$-plan ou un $(\a)$-plan respectivement. Elle n'est donc
pas dans $\R_d$.
Les vari\'et\'es $\R_d^4(a)$ ont pour ouvert dense la famille des
$\pi_L(\Mor{(a,d-a,d)}{\widetilde{\cC}_L})$ au dessus de
$\G$. Celle-ci est irr\'eductible de dimension $3d+7$.
Les vari\'et\'es $\R_d^4(a)$ pour $1\leq a\leq d-1$ sont toutes
distinctes, l'image d'un morphisme g\'en\'eral est donn\'e par une
courbe trac\'ee sur un hyperplan tangent \`a $\G$ qui rencontre un
$(\a)$-plan et un $(\b)$-plan (contenus dans cet hyperplan) en $a$
points et en $d-a$ points.

Remarquons que la dualit\'e \'echange les $(\a)$-plans et les
$(\b)$-plans ce qui montre que les composantes $\R_d^4(a)$ et
$\R_d^4(d-a)$ sont  duales l'une de l'autre. Remarquons enfin que les
strates $\R_{d,1}$ et $\R'_{d,1}$ sont $\R_d^4(1)$ et
$\R_d^4(d-1)$.

\vs 0.2 cm

\noi
{\bf La vari\'et\'e $\R_d^3$}

\vs 0.1 cm

Nous avons vu \`a la proposition \ref{fibrephi} qu'un ouvert dense
$\U_3$ de $\R_d^3$ est form\'e des \'el\'ements $f$ de $\R_d^3$ tels
que la courbe $f(\pu)$ est contenue dans un unique espace lin\'eaire
de codimension $2$ tangent \`a $\G$ (cas $i(f)=1$). Il nous suffit
donc de d\'eterminer les composantes irr\'eductibles de l'ouvert
$\U_3$. 

Nous avons alors un morphisme propre $\Upsilon_3:\U_3\fl{\check T}_\G$
qui a $f$ associe l'unique espace lin\'eaire
de codimension $2$ tangent \`a $\G$ et contenant $f(\pu)$ (l'espace
${\check T}_\G$ est isomorphes aux couples $(x,W)$ o\`u $x\in\G$ et
$W$ est une sous-vari\'et\'e lin\'eaire de dimension 3 de $T_x\G$). Il
s'agit maintenant d'\'etudier la fibre de $\Upsilon_3$. Cette fibre au
dessus de $L_0\in\G$ est isomorphe \`a $\Mor{d}{S_{L_0}}$ o\`u
$S_{L_0}$ est le c\^one de $\G$ d\'ecoup\'e par le sous-espace
lin\'eaire de codimension 2.

\def \cclo {\cC_{L_0}}
\def \qlo {Q_{L_0}}
\def \ecl {\widetilde{\cC}_{L_0}}
\def \N {\mathbb{N}}

Notons $\pi:\widetilde{S}_{L_0}\fl S_{L_0}$ l'\'eclatemement de $S_{L_0}$ au sommet du
c\^one. Nous avons une projection $p$ de $\widetilde{S}_{L_0}$ vers la
droite $L_0$.
Il est facile de v\'erifier que $\widetilde{S}_{L_0}$ est la fibration
en droites projectives au dessus de $L_0$ associ\'ee au faisceau
$\opu(-1)\oplus\opu(1)$. Le groupe de Picard de
$\widetilde{S}_{L_0}$ est donc $\Z^2$ (nous prendrons pour base des
cycles le degr\'e sur $L_0$ et le degr\'e relatif pour $p$).

\begin{lem}\hskip -0.15 cm{\bf .}
Pour tout $\gamma=(a,b)\in\N^2$, si $b\geq 2a$, le sch\'ema
$\Mor{\gamma}{\widetilde{S}_{L_0}}$ est irr\'eductible lisse de
dimension $2a+b+2$, sinon il est vide.
\end{lem}

\dm :
Le th\'eor\`eme principal de \cite{P2} nous permet de dire que le
sch\'ema des morphismes $\Mor{a}{\pu}$ est irr\'eductible et lisse de
dimension $2a+1$. Par ailleurs, la proposition 4 de \cite{P2}
nous permet de dire que comme le degr\'e relatif est positif, alors si
$b\geq 2a$, le sch\'ema $\Mor{\gamma}{\widetilde{S}_{L_0}}$ est
irr\'eductible et lisse de dimension $2a+1+b+1$ et qu'il est vide sinon.

\vs 0.4 cm

L'\'eclatement $\pi$ d\'efinit un morphisme
$\pi:\Mor{\gamma}{\widetilde{S}_{L_0}}\fl\Mor{\pi_*\gamma}{L_0}$. De
plus si $f$ est un \'el\'ement de $\Mor{(a,b)}{\widetilde{S}_{L_0}}$
alors $f(\pu)$ rencontre le diviseur exceptionnel en $\frac{1}{2}b-a$
points, le degr\'e de son image dans $S_{L_0}$ est donc
$\frac{1}{2}b+a$. Ainsi $\pi(f)$ est de degr\'e $d$ si et
seulement si $b=2(d-a)$. Nous supposons maintenant $b\geq 2a$. 

\begin{pr}\hskip -0.15 cm{\bf .}
La vari\'et\'e $\Mor{d}{S_{L_0}}$ a $\d+1$ composantes irr\'eductibles
toutes de dimension $2d+2$ d\'ecrites pas les images par $\pi$ des
$\Mor{(a,2(d-a))}{\widetilde{S}_{L_0}}$ pour $0\leq a\leq\d$.
\end{pr}

\dm :
Comme dans le cas de $\R_d^4$, les morphismes de
$\Mor{(a,b)}{\widetilde{S}_{L_0}}$ dans le sch\'ema $\Mor{d}{S_{L_0}}$
sont des immersions d\`es que $\d\geq a$. Par ailleurs
$\Mor{d}{S_{L_0}}$ est recouvert par ces images (il suffit de regarder
la transform\'ee stricte). Comme toutes ces vari\'et\'es sont
irr\'eductibles de m\^eme dimension elles forment les composantes
irr\'eductibles de $\Mor{d}{S_{L_0}}$.

\vs 0.4 cm

Notons $\R_d^3(a)$ le ferm\'e de $\U_3$ des morphismes $f$ qui sont
dans la composante irr\'eductible
$\pi(\Mor{(a,2(d-a))}{\widetilde{S}_{\Upsilon_3(f)}})$ de
$\Upsilon_3^{-1}(\Upsilon_3(f))$.

\begin{cor}\hskip -0.15 cm{\bf .}
Les vari\'et\'es $\R_d^3(a)$ pour $1\leq a\leq \d$ sont les
composantes irr\'eductibles de $\R_d^3$.  Elles sont de dimension
$2d+9$.\label{compr3}
\end{cor}

\dm :
La vari\'et\'e $\R_d^3(0)$ n'est pas dans $\R_d^3$ (mais dans son
adh\'erence) car un tel morphisme $f$ aurait pour image une droite
contenue dans $\G$, la surface serait alors un c\^one (et sa duale
aussi). Les composantes sont \'evidement autoduales.

\subsection{Compactification de $\R_d$}
\label{bord}

Nous \'etudions dans ce paragraphe la fermeture $\overline{\R_d}$ de
$\R_d$ dans $\p\hom{\LV}{S_d}$. 

\begin{pr}\hskip -0.15 cm{\bf .}
Le bord de $\R_d$ est form\'e de vari\'et\'es isomorphes \`a
$\R_{d'}\times\p(\S_{d-d'})$ pour tout $d'$ entier strictement
inf\'erieur \`a $d$ et de
$\overline{\R_{d,0}}\cup\overline{\R'_{d,0}}$.
\end{pr}

\dm :
Nous avons vu au lemme \ref{defphi} qu'un
\'el\'ement de $\R_d$ est donn\'e par une application lin\'eaire
$\psi$ de $\LV$ dans $\S_d$ telle que le morphisme $\widetilde{\psi}
:\LV\oopu\fl\opu(d)$ est surjectif et que la compos\'ee
$$\opu(-d)\stackrel{^t\!\psi}{\fl}\L^2{\check
V}\oopu\stackrel{q}{\fl}\LV\oopu\stackrel{\psi}{\fl}\opu(d)$$
est nulle (condition not\'ee $(*)$).

\def \wtp {\widetilde{\psi}}

Un \'el\'ement $\psi$ est donc dans le bord de $\R_d$ si il
v\'erifie la condition $(*)$ (qui est ferm\'ee) et si la fl\`eche
$\widetilde{\psi}$ associ\'ee sur les faisceaux n'est plus
surjective. L'\'el\'ement $\psi$ \'etant non nul, l'image de
$\widetilde{\psi}$ dans $\opu(d)$ est un faisceau sans torsion de rang
$1$. C'est donc $\opu(d')$ avec $0\leq d'<d$. La donn\'ee de
$\psi$ est \'equivalente \`a celle d'une surjection de
$\LV\oopu\fl\opu(d')$ et d'une fl\`eche de $\opu(d')$ dans
$\opu(d)$. Ceci correspond au choix d'un \'el\'ement de
$\R_{d'}$ (ou de $\R_{d',0}\cup\R'_{d',0}$) et d'un \'el\'ement de
$\p(S_{d-d'})$. Le cas de $\R_{d',0}\cup\R'_{d',0}$ est adh\'erent
\`a $\R_{d,0}\cup\R'_{d,0}$ d'o\`u le r\'esultat.

\begin{rem}\hskip -0.15 cm{\bf .}
{\rm (\i) Nous pouvons maintenant compl\'eter la preuve du fait
(lemme \ref{defphi}) que la restriction du morphisme $\Phi$ \`a
$\overline{\R_d}$ est d\'efinie en dehors de
$\overline{\R_{d,0}}\cup\overline{\R'_{d,0}}$. En effet, nous avons vu
que $\Phi$ n'est pas d\'efini sur ce ferm\'e et sur le
compl\'ementaire c'est \`a dire sur $\R_d$ et sur les composantes
$\R_{d'}\times\p(S_{d-d'})$ du bord le morphisme $\Phi$ est
d\'efini (voir lemme \ref{phibord} pour l'image du bord).

(\i\i) Notons $i(d,d')$ l'inclusion de $\R_{d'}\times\p(S_{d-d'})$ dans
$\overline{\R_d}$. 
Soient $d'$ et $d''$ deux entiers tels que $0\leq d''<d'<d$. La
compos\'ee de $i(d',d'')\times{\rm{Id}}_{\p(S_{d-d'})}$ et de
$i(d,d')$ qui va de $\R_{d''}\times\p(S_{d'-d''})\times\p(S_{d-d'})$
dans $\overline{\R_{d'}}\times\p(S_{d-d'})$ puis dans
$\overline{\R_d}$ est $i(d,d'')$ o\`u on a identifi\'e
$\p(S_{d'-d''})\times\p(S_{d-d'})$ avec $\p(S_{d-d''})$ gr\^ace \`a la
multiplication $S_{d'-d''}\ot S_{d-d'}\fl S_{d-d''}$.}\label{findefphi}
\end{rem}

\begin{cor}\hskip -0.15 cm{\bf .}
Le bord de $\R_d$ a trois composantes irr\'eductibles : les deux
ferm\'es $\overline{\R_{d,0}}$ et $\overline{\R_{d,0}}$ et la
vari\'et\'e $\overline{\R_{d-1}}\times\p(S_1)$.\label{compbord}
\end{cor}
	
\begin{lem}\hskip -0.15 cm{\bf .}
L'image de $\R_{d'}\times\p(\S_{d-d'})$ par $\Phi$ est la vari\'et\'e
des quadriques contenant la courbe rationnelle normale $C_d$ et dont
le noyau (en tant que forme quadratique) rencontre $C_d$ selon les
$d-d'$ points d\'efinis par l'\'el\'ement de
$\p(\S_{d-d'})$.\label{phibord}
\end{lem}

\dm :
Soient $d-d'$ points de $C_d$ d\'efinis par un \'el\'ement $\xi$ de
$\p(\S_{d-d'})$. L'espace lin\'eaire $H_{\xi}$ engendr\'e par ces points est
donn\'e par le conoyau de $\S_{d'}\stackrel{\xi}{\fl}\S_d$. Soit
$\psi$ un morphisme de $\LV$ dans $S_d$, le noyau de la forme
quadratique $\Phi(\psi)$ contient le sous-espace lin\'eaire $N$
donn\'e par le conoyau de $\psi$. 

Si le morphisme $\psi$ se factorise par un \'el\'ement
$\xi\in\p(\S_{d-d'})$, alors on a une surjection $N\fl H_{\xi}$ ce qui
impose que le noyau de la fl\`eche contient le sous-espace lin\'eaire
$H_{\xi}$ et donc que le noyau de $\Phi(\psi)$ rencontre $C_d$ selon
les $d-d'$ points d\'efinis par $\xi$.

R\'eciproquement, si on a une forme quadratique $w$ qui contient $C_d$
et dont le noyau
rencontre $C_d$ selon les $d-d'$ points d\'efinis par
$\xi\in\p(\S_{d-d'})$, alors l'image de $C_d$ par la projection par
rapport \`a $H_{\xi}$ est la courbe rationnelle normale $C_{d'}$ de
$\p(\S_{d'})$. Son image par la projection par rapport \`a
${\rm{Ker}}w/H_{\xi}$ est alors une courbe rationnelle de degr\'e $d'$ de
$S_d/{\rm{Ker}}w$. En prenant une isom\'etrie de $S_d/{\rm{Ker}}w$
dans $\LV$, on a un \'el\'ement de $\R_{d'}$. Ainsi $w$ est dans
l'image de $\R_{d'}\times\p(\S_{d-d'})$.

\subsection{Degr\'e des images}

Nous \'etudions maintenant les images par $\P$ des vari\'et\'es
$\R_d^k$ ainsi que celle du bord de $\R_d$. Nous donnons notamment
le degr\'e de ces vari\'et\'es et nous d\'ecrirons la courbe
g\'en\'erale de l'image des composantes de $\R_d^4$, $\R_d^3$ et du
bord au prochain paragraphe. 

Pour calculer le degr\'e de l'image, nous utilisons le th\'eor\`eme
\ref{egalite}. Nous savons (proposition \ref{defphi}) que l'image par
$\Phi$ et donc par $\P$ de $\R_d$ est $\aaq_6$. Ceci impose que
l'image de $\R_d^k$ est exactement $\q_k$. Nous v\'erifions que
pour $3\leq k\leq 6$ les vari\'et\'es d\'eterminantielles $\q_k$ ont
les dimensions ``attendues''. Nous pouvons ainsi gr\^ace aux resultats
de \cite{HT} calculer leur degr\'e.

\begin{pr}\hskip -0.15 cm{\bf .}
Les vari\'et\'es $\q_k$ pour $3\leq
k\leq 6$ sont \'equidimensionnelles de dimensions respectives
$(k-2)d-1-\frac{1}{2}(k-1)(k-2)$.\label{equidim}
\end{pr}

\dm :
Nous avons montr\'e que les vari\'et\'es $\R_d^3$, $\R_d^4$, $\R_d^5$
et $\R_d^6$ sont \'equidimension\-nelles de dimensions $2d+9$, $3d+7$,
$3d+8$ et $4d+4$. Les fibres au dessus de $\q_3$, $\q_4$, $\q_5$ et
$\q_6$ \'etant \'equidimensionnelles de dimensions $d+11$, $d+11$,
$15$ et $15$ (proposition \ref{fibrephi}), les vari\'et\'es $\q_k$
pour $3\leq k\leq 6$ sont donc \'equidimensionnelles de la dimension
voulue.

\vs 0.2 cm

Nous pouvons maintenant montrer le th\'eor\`eme suivant :

\begin{theo}\hskip -0.15 cm{\bf .}
L'image par $\P$ de la vari\'et\'e $\R_d$ (resp. de $\R_d^5$) est un
localement ferm\'e irr\'eductible de dimension $4d-11$ (resp. $3d-7$)
et de degr\'e $i(d)$ (resp. $j(d)$) qui est $PGL_2$-invariant.

L'image par $\P$ de la vari\'et\'e $\R^4_d$ (resp. $\R^3_d$) est de
degr\'e $k(d)$ (resp. $p(d)$). Elle  a $\d$ composantes
irr\'eductibles qui sont des localement ferm\'es (resp. ferm\'es) de
dimension $2d-4$ (resp. $d-2$) invariants sous $PGL_2$.

Si $d\geq 5$, l'image du bord de $\R_d$ priv\'e de $\overline{\R_{d,0}}$ et
$\overline{\R'_{d,0}}$ est irr\'eductible de dimension $4d-14$ et de
degr\'e $2(4d-14)i(d-1)$.

Les degr\'es sont donn\'es par :
\vs -0.6 cm
$$i(d)=\prod_{k=0}^{d-6}\frac{\binom{d+1+k}{d-5-k}}{\binom{2k+1}{k}} \
\ \  
j(d)=\prod_{k=0}^{d-5}\frac{\binom{d+1+k}{d-4-k}}{\binom{2k+1}{k}}
 \ \
\  
k(d)=\prod_{k=0}^{d-4}\frac{\binom{d+1+k}{d-3-k}}{\binom{2k+1}{k}}
\ \
\ 
p(d)=\prod_{k=0}^{d-3}\frac{\binom{d+1+k}{d-2-k}}{\binom{2k+1}{k}}$$
\label{degre}
\end{theo}

\vs -0.3 cm

\dm :
Les r\'esultats de J. Harris et L.W. Tu \cite{HT}, que l'on peut
appliquer aux vari\'et\'es d\'eterminantielles
$\q_k$, nous permettent de calculer les
degr\'es.

Le r\'esultat sur $\R_d$ et $\R_d^5$ d\'ecoule directement de la
proposition \ref{irredrd} et de la dimension des fibres (proposition
\ref{fibrephi}).

Les $d-1$ composantes irr\'eductibles $\R_d^4(a)$ pour $1\leq a\leq
d-1$ de $\R_d^4$ sont telles que $\R_d^4(a)$ et $\R_d^4(d-a)$ sont en
dualit\'e (corollaire \ref{compr4}). Leurs images sont donc les
m\^emes ce qui donne $\d$ ferm\'es irr\'eductibles de
$\P(\R_d^4)$. Nous montrerons \`a la proposition \ref{poncelet} que
ces ferm\'es sont distincts et nous les d\'ecrirons.

Nous avons vu que l'image de $\R_d^3$ est $\q_3$ qui est ferm\'e (car $\q_2$
est vide) de dimension $d-2$ (proposition \ref{equidim}). Chacune des
composantes irr\'eductibles de $\R_d^3$ donnera un ferm\'e irr\'eductible de
$\q_3$. Nous montrerons \`a la proposition \ref{poncelet} que ces
ferm\'es sont distincts et nous les d\'ecrirons.

Le bord de $\R_d$ priv\'e de $\overline{\R_{d,0}}$ et
$\overline{\R'_{d,0}}$ est $\overline{\R_{d-1}}\times\p(\S_1)$
(corollaire \ref{compbord}). Il est donc irr\'eductible de dimension
$4d+1$. D'apr\`es la proposition \ref{fibrephi}, son image est de
dimension $4d-14$. Nous d\'ecrirons son image au lemme \ref{imbordpsi}
et nous donnerons son degr\'e.

\subsection{Description g\'eom\'etrique de quelques images}

Nous d\'eterminons g\'eom\'etriquement les courbes de l'image du bord,
de $\P(\R_d^4)$ et de $\P(\R_d^3)$.

\begin{rem}\hskip -0.15 cm{\bf .}
{\rm Nous pouvons d\'ecrire g\'eom\'etriquement la correspondance entre les
formes quadratiques de $\S_d$ contenant $C_d$ et les courbes de
degr\'e $d-2$ de $\p(\S_2)$. Soit $Q$ une forme quadratique de $S_d$
contenant $C_d$ et d\'efinissons la courbe de $S^2C_d$ suivante : 
\vs -0.5 cm
$$C(Q)=\{(p,q)\in S^2C_d\  /\ {\rm{la}}\ {\rm{droite}}\  (pq)\
{\rm{est}}\ {\rm{isotrope}}\ {\rm{pour}}\ Q\}$$ 
\vs -0.1 cm
La vari\'et\'e $S^2C_d$ est canoniquement isomorphe \`a $S^2\pu$ et le
morphisme :
\vs -0.5 cm
$$C:\p(S^2S_{d-2})\fl\p(S^{d-2}S_2)$$
\vs -0.9 cm
$$Q\mapsto C(Q)$$
est la description g\'eom\'etrique de la loi de r\'eciprocit\'e de
Hermite. Ainsi si nous consid\'erons $\Phi$ comme \`a valeur dans
$\p(S^2S_{d-2})$, nous avons $\P=C\circ\Phi$. En effet, soit $\psi$
une projection de $\p(S_d)$ dans $\plv$. Les droites $L_p$ et $L_q$
correspondant aux images de $p$ et $q$ se coupent si et seulement si
la droite $(\psi(p)\psi(q))$ est contenue dans $\G$ ce qui signifie
exactement que $(pq)$ est isotrope pour $Q$.}
\end{rem}

\begin{lem}\hskip -0.15 cm{\bf .}
L'image de $\R_{d'}\times\p(\S_{d-d'})$ est form\'ee des courbes
r\'eunion d'une courbe de degr\'e $d'-2$ de $\P(\R_{d'})$ et des
$d-d'$ tangentes \`a la conique $C_0$ d\'efinies par l'\'el\'ement
de $\p(\S_{d-d'})$. Elle est de degr\'e
$2^{d-d'}\binom{d+3d'-11}{d-d'}i(d')$ pour $d\geq 5$.\label{imbordpsi}
\end{lem}

\dm :
Si $\psi$ est dans le bord de $\R_d$, disons que
$\psi=(\psi',\xi)\in\R_{d'}\times\p(S_{d-d'})$, alors $\psi$ se
factorise par $\xi$ ce qui impose (voir le lemme \ref{phibord}) que le
noyau de la forme quadratique $Q=\Phi(\psi)$ rencontre $C_d$ selon les
$d-d'$ points d\'efinis par $\xi$. Si $p$ est un tel point, alors pour
tout $q\in C_d$ la droite $(pq)$ est isotrope pour $Q$. Ainsi les
paires $(p,q)$ pour tout $q\in C_d$ sont dans $C(Q)$. Le lieu
singulier abstrait contient donc les tangentes \`a la conique
canonique d\'efinies par $\xi$. La composante restante est d\'efinie
par $\psi'$ et est un \'el\'ement de $\P(\R_{d'})$.

R\'eciproquement, si la courbe $C(Q)$ contient les droites tangentes
\`a $C_0$ d\'efinies par un \'el\'ement
$\xi\in\p(S_{d-d'})$, alors pour tout point $p$ d\'efini par $\xi$,
nous savons d'apr\`es la remarque pr\'ec\'edente que pour tout $q\in
C_d$ la droite $(pq)$ est isotrope pour $Q$. Ceci est \'equivalent \`a
dire que $p$ est dans le noyau de $Q$. Il suffit alors de prendre une
isom\'etrie de $S_d/{\rm{KerQ}}$ dans $\LV$ pour trouver un
\'el\'ement de $\R_{d'}\times\p(S_{d-d'})$ qui s'envoie sur la courbe
$C(Q)$.

La fibre des courbes contenant une droite tangente n'est pas contenue
dans le bord de $\R_d$. Par exemple si $f$ est un \'el\'ement de
$\U_4$ tel que $f(\pu)$ rencontre le sommet du c\^one
$\cC_{\Upsilon_4(f)}$ alors son lieu singulier abstrait contient une
droite tangente \`a $C_0$.

Pour calculer son degr\'e, il suffit de calculer le nombre de courbes
de ce ferm\'e qui passent par $d+3d'-11$ points en position
g\'en\'erale. Mais on connait le degr\'e de $\P(\R_{d'})$ dans la
vari\'et\'e des courbes de degr\'e $d'-2$, c'est $i(d')$ (il
repr\'esente le nombre de courbes de ce type passant par $4d'-11$
points). Pour d\'eterminer une courbe passant par $3d'+d-11$ points il
suffit d'en choisir $d-d'$ ($\binom{d+3d'-11}{d-d'}$ possibilit\'es)
qui nous donnent les tangentes ($2^{d-d'}$ possibilit\'es) et les
$4d'-11$ autres nous d\'efinissent $i(d')$ courbes de degr\'e $d'-2$
dans l'ensemble voulu. Le degr\'e recherch\'e est donc
$2^{d-d'}\binom{d+3d'-11}{d-d'}i(d')$.

Si $d=4$, la vari\'et\'e $\R_3$ ne v\'erifie plus les conditions
pr\'ec\'edentes. La dimension du bord de $\R_4$ est $18$, celle de son
image est $3$ et son degr\'e est $6$.

\begin{lem}\hskip -0.15 cm{\bf .}
L'image par $\P$ de $\R_d^4(a)$ est contenue dans la vari\'et\'e des
courbes r\'eunion d'une courbe de $\crp_{a-1}$ et d'une courbe de
$\crp_{d-a-1}$.\label{degponc}
\end{lem}

\dm :
Il suffit de v\'erifier ce lemme sur l'ouvert $\U_4$ de $\R_d^4$ et
m\^eme sur l'ouvert de $\R_d^4(a)$ des courbes ne passant pas par le
sommet du c\^one. Si $f\in\R_d^4(a)$, le sommet du c\^one
$\Upsilon(f)$ est une droite de $\pv$. Notons $W$ le quotient de
dimension $2$ de $V$ qui la d\'efinit et $N_W$ le sous-espace de
dimension $2$ de $V$ correpondant. 

Le quotient ${\cal L}$ de la fl\`eche $N_W\oopu\fl f^*Q$ est
localement libre de rang $1$ si $f(\pu)$ ne passe pas par le sommet
(voir preuve de la proposition \ref{adhe}). Par ailleurs, ce faisceau
d\'ecrit la trace des g\'en\'eratrices sur la droite $\p(W)$. Son
degr\'e est donc donn\'e par le degr\'e de $f(\pu)$ sur le premier
facteur de la quadrique $Q_{\p(W)}$ qui s'identifie \`a
$\p(W)\times\p({\check N_W})$. Ainsi ${\cal L}$ est isomorphe \`a
$\opu(a)$. Nous avons donc un quotient naturel $S_a$ de $H^0(f^*Q)$
(dont le noyau est $S_{d-a}$) tel que l'application lin\'eaire $V\fl
H^0(f^*Q)\fl S_a$ se factorise par $W$.

Soient $N$ le noyau et $Q$ le conoyau de la fl\`eche de $V\oops$ dans
$K_a$ (voir notations) d\'eduite de la compos\'ee pr\'ec\'edente, le
faisceau $R$ de support le lieu singulier abstrait est donn\'e par la
suite exacte :
$$0\fl N \fl K_{d-a}\fl R\fl Q\fl 0$$
Mais l'image de $V\oops$ dans $K_a$ est $W\oops$ donc $Q$ est donn\'e
par :
$$0\fl W\oops\fl K_a\fl Q\fl 0$$
et $N=N_W\oops$. Ainsi $Q$ a pour support une courbe de
$\crp_{a-1}$. Enfin, le noyau de l'application de $R$ dans $Q$ est
donn\'e par le conoyau de la fl\`eche $N_W\oops \fl K_{d-a}$ ce qui
nous dit que son support est une courbe de $\crp_{d-a-1}$. 

Si l'\'el\'ement $f$ de $\R_d^4(a)$ est tel $f(\pu)$ passe par le
sommet du c\^one, alors le conoyau de la fl\`eche $N_a\oopu\fl f^*Q$
est un faisceau ayant de la torsion et la courbe de $\crp_{d-a-1}$ se
d\'ecompose en une courbe de $\crp_{d'-a-1}$ et les $d-d'$ droites
tangentes \`a la conique d\'etermin\'ees par les points de torsion de
ce faisceau.

\begin{lem}\hskip -0.15 cm{\bf .}
La vari\'et\'e $\crp_n$ est un ferm\'e irr\'eductible et r\'eduit de
dimension $2n$ et de degr\'e $l(n)$ de $\p(S^n\S_2)$.
\end{lem}

\dm :
Les courbes de degr\'e $n$ en relation de Poncelet avec la conique
canonique sont donn\'ees par les pinceaux de sections de
$H^0K_{n+1}$. Elles s'identifient ainsi \`a
${\rm{Grass}}(2,\S_{n+1})$ dans $\p(\L^2\S_{n+1})=\p(S^n\S_2)$ et
forment donc une vari\'et\'e irr\'eductible et r\'eduite de dimension
$2n$ et de degr\'e :
\vs -0.5 cm
$$l(n)=\frac{\binom{2n+1}{n+1}}{2n+1}=\frac{(2n)!}{n!(n+1)!}\ \ \ \ .$$

\begin{pr}\hskip -0.15 cm{\bf .}
Les vari\'et\'es $\R_d^4(a)$ et $\R_{d}^4(d-a)$ ont la m\^eme image
par $\P$, elle est irr\'eductible, r\'eduite et isomorphe \`a
$\crp_{a-1}\times\crp_{d-a-1}$. Son degr\'e est $m(a,d)$ que l'on calculera.

L'image r\'eduite de $\R_d^3(a)$ est isomorphe \`a
$\crp_{a-1}\times\p(\S_{d-2a})$ et est form\'ee des courbes r\'eunion
de $d-2a$ tangentes \`a $C_0$ et d'une courbe de
$\crp_{a-1}$ compt\'ee deux fois. Elle est de degr\'e
$2^{d-2a}l(a-1)$.\label{poncelet}
\end{pr}

\dm :
Les vari\'et\'es $\R^4_d(a)$ et $\R^4_{d}(d-a)$ \'etant
irr\'eductibles r\'eduites et \'echan\-g\'ees par dualit\'e, elles ont
la m\^eme image par $\P$ qui est irr\'eductible et r\'eduite. De plus,
cette image est de dimension $2d-4$. On a vu au lemme \ref{degponc} que cette
image est contenue dans $\crp_{a-1}\times\crp_{d-a-1}$. Le lemme
pr\'ec\'edent nous permet de conclure \`a l'\'egalit\'e.

On peut calculer le degr\'e $m(a,d)$ de l'image de $\R_d^4(a)$ : il
suffit de calculer le nombre de courbes r\'eunion d'une courbe de
$\crp_{a-1}$ et d'une courbe de $\crp_{d-a-1}$ passant par $2d-4$
points. On choisit $2a-2$ points ($(_{2a-2}^{2d-4})$ choix) et par ces
points il passe $l(a-1)$ courbes de $\crp_{a-1}$, par les $2d-2a-2$
points restant il passe $l(d-a-1)$ courbes de $\crp_{d-a-1}$. Le
degr\'e est donc $m(a,d)=(_{2a-2}^{2d-4})l(a-1)l(d-a-1)$ sauf dans le
cas particulier o\`u $a=\frac{d}{2}$ car on a compt\'e deux fois
chaque situation et $m(\frac{d}{2},d)=\frac{1}{2}(_{\
d-2}^{2d-4})l(\frac{d}{2}-1)^2$. 

La construction du lemme \ref{degponc} est encore valable pour une
courbe de $\R_d^3(a)$. Le conoyau de la fl\`eche de $W\ot\ops$ dans
$K_a$ est le m\^eme et donne une courbe $X\in\crp_{a-1}$. La fl\`eche
de $N_W\ot\ops$ dans $K_{d-a}$ se d\'ecompose en ${\check W}\ot\ops\fl
K_a\fl K_{d-a}$ ce qui nous donne la courbe $X$ une deuxi\`eme fois
(conoyau de ${\check W}\ot\ops\fl K_a$) et $d-2a$ droites tangentes
\`a $C_0$ (conoyau de $K_a\fl K_{d-a}$). On conclue par
irr\'eductibilit\'e car l'image de $\R_d^3(a)$ et la vari\'et\'e
$\crp_{a-1}\times\p(\S_{d-2a})$ ont la m\^eme dimension.

\begin{rem}\hskip -0.15 cm{\bf .}
{\rm Une surface g\'en\'erale a pour lieu singulier
abstrait une courbe non singuli\`ere. En effet, l'image de $\R_d^4(1)$
est $\crp_{d-2}$. Or il existe des courbes lisses dans $\crp_{d-2}$,
par g\'en\'eralisation on a le r\'esultat.}
\end{rem}

\begin{pr}\hskip -0.15 cm{\bf .}
Une surface est dans $\R_d^4(a)$ si et seulement si son lieu singulier
est r\'eunion d'une droite $a$-uple et d'une courbe de degr\'e
$\frac{1}{2}(d-a-1)(d+a-2)$ et celui de sa duale est r\'eunion d'une
droite $(d-a)$-uple et d'une courbe de degr\'e
$\frac{1}{2}(a-1)(2d-a-2)$.
\end{pr}

\dm :
Avec les notations ci-dessus la droite $\p({\check N_W})$ est multiple
pour la duale. La multiplicit\'e de cette droite est le nombre de
g\'en\'eratrices de la surface contenues dans le plan correspondant
\`a un point de $\p({\check N_W})$. Ceci revient \`a calculer, pour un
sous-espace vectoriel ${\bC}$ de dimension $1$ de $N_W$ donn\'e le
nombre de points de $\pu$ tels que la fl\`eche induite ${\bC}\oopu\fl
f^*Q$ soit nulle. Si $f$ est dans $\R_d^4(a)$, cette fl\`eche se
factorise par $\opu(d-a)$ (car $N_W\oopu$ se factorise par ce
faisceau) il suffit donc qu'elle s'annule dans $\opu(d-a)$. On a donc
$d-a$ tels points. Ainsi la duale a une droite $(d-a)$-uple.

De plus la duale ${\check S}$ d'une surface $S$ de $\R_d^4(a)$ est
dans $\R_d^4(d-a)$. Nous savons donc que ${({\check S}\check )}$ a une
droite $a$-uple. Ainsi $S$ a une droite $a$-uple. Pour calculer le
degr\'e de la courbe r\'esiduelle, il suffit d'utiliser les
r\'esultats de \cite{KL}.

Si le lieu singulier abstrait contient une courbe de $\crp_n$ alors
les $n+1$ points de la conique d\'efinissant un polyg\^one associ\'e
correspondent \`a $n+1$ g\'en\'eratrices de la surface qui se coupent
deux \`a deux. Ceci n'est possible que si elles sont concourrantes ou
coplanaires. On a alors une courbe $(n+1)$-uple sur la surface ou sa
duale. Si on a une droite $(n+1)$-uple sur la
surface ou sa duale, le lieu singulier abstrait contient alors une
courbe de $\crp_n$.

R\'eciproquement, si on a une surface dont le lieu singulier v\'erifie les
hypoth\`eses de la proposition, alors son lieu singulier abstrait
contient une courbe de $\crp_{a-1}$ (droite $a$-uple sur la surface)
et une courbe de $\crp_{d-a-1}$ (droite $(d-a)$-uple de la
duale). Pour une raison de degr\'e, le lieu singulier abstrait est
r\'eunion de ces deux courbes. Le lieu singulier de cette surface est
donc (proposition
\ref{poncelet}) contenue dans $\P(\R_d^4(a))$. Cependant, nous savons
(proposition \ref{fibrephi}) que la fibre au dessus d'un point de
$\P(\R_d^4(a))$ a deux composantes irr\'eductibles, l'une contenue dans
$\R_d^4(a)$ l'autre dans $\R_d^4(d-a)$. Comme la surface contient une
droite $a$-uple elle est dans $\R_d^4(a)$.


\begin{rem}\hskip -0.15 cm{\bf .}
{\rm Si $f\in\R_d^4(a)\cap\R_d^4(d-a)$ alors $f(\pu)$ est trac\'ee sur une
quadrique $L_1\times L_2$ de $\G$ d\'efinie par deux droites $L_1$ et
$L_2$ de $\pv$ (c'est l'ensemble des droites qui rencontrent les deux
droites). La courbe $f(\pu)$ est de bidegr\'e $(a,d-a)$ sur la
quadrique. Elle a donc $(a-1)(d-a-1)$ points doubles.
La surface est alors autoduale et son lieu singulier est r\'eunion de la
droite $a$-uple $L_1$ de la droite $(d-a)$-uple $L_2$ et des
$(a-1)(d-a-1)$ g\'en\'eratrices correspondant aux points
doubles.}
\end{rem}


\subsection{\'Etude des surfaces de degr\'e $5$}

Nous d\'ecrivons ici les fibres du morphisme
${\cal M}(d)\stackrel{\P}{\fl}\p(S^{d-2}S_2)$ de fa\c con
g\'eom\'etrique dans le
cas des surfaces de degr\'e $5$. Ceci nous permet de retrouver des
r\'esultats de \cite{P1} sur la position relative d'une cubique et
d'une conique.

%

Notons ${\cal M}^s(5)$ l'image dans ${\cal M}(5)$ de l'ouvert des
points stable, $\P({\cal M}^s(5))$ s'identifie \`a $\q_6\cup\q_5$.

\begin{pr}\hskip -0.15 cm{\bf .}
La fibre de l'application ${\cal M}^s(5)\stackrel{\P}{\fl}\q_6\cup\q_5$ au
dessus d'une cubique est donn\'ee par les triangles de Poncelet
associ\'es \`a la cubique.\label{fibredeg5}
\end{pr}

\dm :
Les surfaces de ${\cal M}^s(5)$ sont de type $2$ : les surfaces de
type $1$ ou dont la duale est de type $1$ sont toutes dans $\R_5^4$
(voir corollaire \ref{compr4}). Le faisceau $R$ d\'efinissant le lieu
singulier abstrait est alors le conoyau du morphisme injectif
$V\oops\fl K_2\oplus K_3$. De plus, le morphisme canonique (invariant
sous le groupe $G_2$) de $V$ dans $H^0K_2$ est n\'ecessairement
surjectif. Ainsi on a la suite exacte :
\vs -0.5 cm
$$0\fl N\fl K_3\fl R\fl 0\ \ \ (*)$$
o\`u le faisceau $N$ est $\ops(-1)\oplus\ops$ (le faisceau $V\oops$
est une extension non triviale de $K_2$ par $N$). Ceci nous d\'efinit
une section de $K_3$ (c'est l'intersection $V\cap H^0K_3$ dans
$K_2\oplus K_3$) qui nous permet d'\'ecrire la suite exacte :
\vs -0.4 cm
$$0\fl \ops(-1)\stackrel{C}{\fl}\I_Z(2)\fl R\fl 0 \ \ \ (**)$$
o\`u $Z$ est la r\'eunion des sommets d'un triangle dont les c\^ot\'es
sont tangents \`a $C_0$ et $C$ est l'\'equation de la
cubique support de $R$.

Ainsi, \`a toute surface $S$ de ${\cal M}^s(5)$ on associe la cubique
$\P(S)$ et le triangle de Poncelet associ\'e \`a $\P(S)$ d\'efini par
la section $V\cap H^0K_3$ de $K_3$ (qui est invariante sous l'action
de $G_2$).

R\'eciproquement, si on a une cubique $C$ et un triangle de Poncelet
associ\'e de sommets $Z$, alors la suite exacte $(**)$ permet de
retrouver le faisceau $R$. Pour retrouver le sous-espace vectoriel $V$
de $H^0(K_2\oplus K_3)$, il suffit de d\'efinir un morphisme de
$K_2\oplus K_3$ dans $R$ ($V$ sera le noyau sur les sections). De tels
morphismes prolongeant celui de $K_3$ dans $R$ d\'efini par $C$ et $Z$
(suite exacte $(*)$) sont param\'etris\'es par
${\rm{Hom}}(K_2,R)$. Cet espace vectoriel est extension de
${\rm{Ext}}^1(K_2,\ops(-1))=S_0$ par ${\rm{Hom}}(K_2,K_3)=\S_1$. Le
sous-groupe unipotent de $G_2$ est isomorphe \`a $S_1$ et agit
transitivement sur $\hom{K_2}{K_3}$. Le quotient $G_2/S_1$ qui est
alors isomorphe \`a $\bC^*$ agit par homoth\'eties sur
${\rm{Ext}}^1(K_2,N)$. Ainsi le groupe $G_2$ a deux orbites dans
${\rm{Hom}}(K_2,R)$ dont l'une est le vecteur nul. Comme $V\oops$ est
une extension non triviale de $K_2$ par $N$, cette orbite ne peut
convenir. L'autre orbite sous $G_2$ nous d\'efini un \'el\'ement de la
fibre de $C$. 



\begin{cor}\hskip -0.15 cm{\bf .}
Il existe exactement deux triangles de Poncelet associ\'es \`a une
cubique g\'en\'erale de $\ps$. Il existe une hypersurface de degr\'e
$6$ de $\p(S^3S_2)$ telle que la cubique g\'en\'erale a un unique
triangle de Poncelet associ\'e. Il existe trois ferm\'es
irr\'eductibles de codimension $3$ et de degr\'es $5$, $30$ et $12$ de
$\p(S^3S_2)$ tels que la cubique g\'en\'erale a une infinit\'e de
triangles associ\'es.
\end{cor}

\dm :
On sait que la fibre de l'application de ${\cal M}(5)$ dans
$\p(S^3\S_2)$ est donn\'ee par deux points au dessus de $\q_6$, un
seul au dessus de $\q_5$ et une infinit\'e de points au dessus de
$\q_4$. La proposition pr\'ec\'edente permet de trouver les deux
premiers cas. De plus l'image de $\R_5^4$ est form\'ee de deux
composantes : $\crp_3$ et les courbes r\'eunion d'un \'el\'ement de
$\crp_2$ et d'une droite. Ces deux types de courbes ont une infinit\'e
de triangles associ\'es et forment deux ferm\'es de codimension $3$ de
degr\'es $5$ et $30$. Enfin l'image du bord de $\R_5$ est form\'e des
courbes r\'eunion d'une conique et d'une droite tangente \`a $C_0$. Il y a une infinit\'e de triangles associ\'es \`a
ces courbes. Cette vari\'et\'e est de degr\'e $12$.

\subsection{Une param\'etrisation birationnelle de l'image de $\P$}

Dans ce paragraphe, nous donnons une param\'etrisation
birationnelle de l'image de $\P$ qui permet de montrer que l'image est
rationnelle. Nous d\'ecrivons explicitement cette param\'etrisation.

\begin{pr}\hskip -0.15 cm{\bf .}
Si $d$ est pair, l'image de $\P$ est birationnelle \`a la vari\'et\'e
des formes quadratiques de rang $4$ de $\S_{d-2}$.\label{birat}
\end{pr}

\dm :
Supposons que $d$ vaut $2n$. Consid\'erons l'ouvert
$\R_{2n,n}\cap\R'_{2n,n}$. Si $f$ est un
\'el\'ement de cet ouvert alors $f^*K$ et $f^*Q$ s'identifient \`a
$U\ot\opu(-n)$ et $W\oopu(n)$ (o\`u $U$ et $W$ sont des espaces
vectoriels de dimension $2$). La condition de parit\'e est
indispensable pour que $f^*K$ et $f^*Q$ aient ces d\'ecompositions.

Fixons $U$ et $W$, nous allons montrer que l'image de $\P$ est
birationnelle au bon quotient de $\p\hom{W\ot {\check U}}{\S_{d-2}}$
par $PGL(U)\times PGL(W)$. Cet espace projectif est isomorphe \`a
l'espace $\p{\rm {Ext}}^1(W\ot\opu(n),U\ot\opu(-n))$. Ses points
semi-stable pour l'action de $PGL(U)\times PGL(W)$ sont d\'efinis par la
surjectivit\'e des morphismes $W\otimes\S_{d-2}\fl U$ et ${\check
U}\otimes\S_{d-2}\fl {\check W}$.


Un \'el\'ement de $\p{\rm {Ext}}^1(W\ot\opu(n),U\ot\opu(-n))$
d\'efinit un faisceau $F$ de rang $4$ sur $\pu$. Le faisceau $F$
est trivial si et seulement si la fl\`eche $W\ot\S_{n-1}\fl
U\ot\S_{n-1}$ est un isomorphisme (ce qui d\'ecrit un ouvert plus petit
que celui des points semi-stables). Nous pouvons appliquer la
construction de la d\'efinition \ref{deflieu} sur cet ouvert pour
construire une courbe $X$ de degr\'e $d-2$ de $\ps$. En identifiant
$F$ \`a $V\oopu$, nous constatons que l'\'el\'ement de
$\R_{2n,n}\cap\R'_{2n,n}$ ainsi d\'efini a $X$ pour image par
$\P$. Nous avons donc construit un morphisme qui a une extension
triviale associe la courbe $X$ (voir remarque \ref{calcexp} pour une
description explicite). Elle est \'evidement invariante sous
$PGL(U)\times PGL(W)$, et nous donne une application rationnelle $\pi$
du bon quotient de $\p\hom{U\ot {\check W}}{\S_{d-2}}$ par
$PGL(U)\times PGL(W)$ dans l'image de $\P$.


R\'eciproquement, la donn\'ee de $f$ dans $\R_{2n,n}\cap\R'_{2n,n}$
d\'efinit la classe modulo l'action de $PGL(U)\times PGL(W)$ de
l'extension suivante :
\vs -0.5 cm
$$0\fl U\ot\opu(-n)\fl V\oopu\fl W\oopu(n)\fl 0$$ 
qui est un \'el\'ement de $\p{\rm
{Ext}}^1(W\ot\opu(n),U\ot\opu(-n))$.
Le fait que le terme central est trivial impose que la fl\`eche
$W\ot\S_{n-1}\fl
U\ot\S_{n-1}$ d\'eduite de l'extension est un isomorphisme. On a donc
une extension semi-stable. L'application qui a $f$ associe la classe
de cette extension modulo $PGL(U)\times PGL(W)$ est donc \`a valeur
dans le bon quotient. De plus, elle est \'evidement invariante sous
$PGL(V)$ et par dualit\'e (identification entre $U$ et ${\check
W}$). Ainsi elle se factorise par l'image de $\P$ (au moins sur
l'ouvert des points stables). Ceci nous montre que l'application
rationnelle $\pi$ est g\'en\'eriquement bijective. La normalit\'e de
$\P(\R_d)$ permet de conclure ces deux vari\'et\'es sont
birationnelles.

Enfin le quotient pr\'ec\'edent est birationnel \`a la vari\'et\'e
des formes quadratiques de rang $4$ de
$S_{d-2}$. En effet, notons $q'$ la forme quadratique de rang $4$ sur $W\ot
{\check U}$ et consid\'erons l'application suivante :
$${\rm{Hom}}(W\ot {\check U},\S_{d-2})\fl S^2\S_{d-2}$$
\vs - 0.8 cm
$$\hskip 1 cm \psi\mapsto\psi {q'}^{t}\!\psi$$
L'image de ce morphisme est la vari\'et\'e des formes quadratiques de rang
au plus $4$ et la fibre g\'en\'erale du morphisme induit sur les
espaces projectifs est $PGL(U)\times PGL(W)$. Le quotient est donc
birationnel \`a la vari\'et\'e des formes quadratiques de rang
$4$.


 
\begin{cor}\hskip -0.15 cm{\bf .}
L'image de $\P$ est rationnelle
\end{cor}

\begin{rem}\hskip -0.15 cm{\bf .}
{\rm Nous pouvons calculer explicitement le morphisme $\pi$ d\'ecrit dans
la proposition pr\'ec\'edente. En effet, soit $\vp$ un \'el\'ement
$\vp$ de $\p\hom{U\ot{\check W}}{S_{d-2}}$. L'image de $\pi$ est
donn\'ee par le d\'eterminant de la compos\'ee suivante :
$W\ot\S_{n-2}\ot\S_2\fl W\ot\S_n\stackrel{\vp}{\fl}U\ot\S_{n-2}$ vue
comme matrice de taille $2(n-1)\times2(n-1)$ \`a coefficients dans
$\S_2$. 
 
Consid\'erons $\vp$ comme une matrice
$(a_{i,j})_{(i,j)\in [1,2]\times[1,2]}$ \`a coefficients
dans $\S_{d-2}$ (les polyn\^omes de degr\'e $d-2$ en deux
variables). Si $Q$ est un \'el\'ement de $S_k$, notons ${\rm d}_Q$
l'application $S_{x+k}\fl S_x$ d\'efinie par $Q$ (qui vient de
$S_{x+k}\ot S_k\fl S_x$). Alors l'image de $\vp$ par $\pi$ est
donn\'ee par le d\'eterminant de la matrice $2(n-1)\times 2(n-1)$
donn\'ee comme une matrice $2\times 2$ dont les blocs sont :
$\big(({\rm d}_{x^{k+l}y^{2n-4-(k+l)}}(a_{i,j}))_{(k,l)\in
[0,n-2]\times[0,n-2]}\big)_{(i,j)\in [1,2]\times[1,2]}$ pour plus
de d\'etails sur les calculs dans les repr\'esentations de $SL_2$ voir
\cite{SP}.}
\label{calcexp}
\end{rem}

\begin{rem}\hskip -0.15 cm{\bf .}
{\rm Consid\'erons l'espace projectif $\p(S^2S_{d-2})$ comme celui des
formes quadratiques sur $S_d$ qui contiennent la courbe rationnelle
normale $C_d$ (cf. remarque \ref{memimage}), soit $Q$ une telle forme
quadratique, supposons que $C_d$ ne rencontre pas le lieu singulier
de $Q$. Notons $Q_0$ le lieu lisse de $Q$, le fibr\'e tangent
$T_{Q_0}(1)$ est muni d'une forme quadratique de rang \'egal \`a ${\rm
rg}(Q)-2$.

Supposons maintenant que $d$ est pair, alors un \'el\'ement
g\'en\'eral $Q$ de $\p(S^2S_{d-2})$ v\'erifie le fait que
$T_{Q_0}(-1)\vert_{C_d}$ est trivial et on a une identification
$S_{d-2}\fl H^0T_{Q_0}(-1)\vert_{C_d}$. Ce n'est plus vrai si $d$ est
impair. Ainsi, dans le cas pair nous avons d\'efini une
application rationnelle $p$ de $\p(S^2\S_{d-2})$ l'espace des formes
quadratiques de $S_d$ contenant la courbe rationnelle normale vers
$\p(S^2\S_{d-2})$ l'espace des formes quadratiques de $S_{d-2}$. De
plus, l'image d'une forme quadratique de rang $k$ est une forme
quadratique de rang $k-2$. 

Remarquons enfin que dans le cas des formes quadratiques de rang $6$,
l'application rationnelle $p$ est la r\'eciproque de l'application
rationnelle $\pi$ de la proposition \ref{birat}. Elle est donc
birationnelle. De m\^eme dans le cas des rang inf\'erieurs \`a 6, les
isomorphismes exceptionnels permettent de montrer que $p$ est
birationnelle. Je ne sais pas si c'est le cas en g\'en\'eral.}
\end{rem}

\end{document}